\documentclass[global]{svjour}
\usepackage{amsfonts}
\usepackage{amsmath}
\usepackage{amssymb}
\usepackage{graphicx}

\input{tcilatex}
\begin{document}

\journalname{}
\title{Approximate Controllability for Linear Stochastic Differential
Equations in Infinite Dimensions}
\author{D. Goreac\inst{1}}
\institute{Laboratoire de Math\'{e}matiques, Unit\'{e} CNRS UMR 6285,\\
Universit\'{e} de Bretagne Occidentale,\\
6, av. Victor LeGorgeu, B.P. 809, \\
29200 Brest cedex, France\\
\email{Dan.Goreac@univ-brest.fr}%
\\
Tel. 02.98.01.72.45, Fax. 02.98.01.67.90\\
\textbf{AMS} \textbf{Classification.} 60H10,60H15}
\maketitle

\begin{abstract}
The objective of the paper is to investigate the approximate controllability
property of a linear stochastic control system with values in a separable
real Hilbert space. In a first step we prove the existence and uniqueness
for the solution of the dual linear backward stochastic differential
equation. This equation has the particularity that in addition to an
unbounded operator acting on the $Y$-component of the solution there is
still another one acting on the $Z$-component. With the help of this dual
equation we then deduce the duality between approximate controllability and
observability. Finally, under the assumption that the unbounded operator
acting on the state process of the forward equation is an infinitesimal
generator of an exponentially stable semigroup, we show that the generalized
Hautus test provides a necessary condition for the approximate
controllability. The paper generalizes former results by Buckdahn,
Quincampoix and Tessitore (2006) and Goreac (2007) from the finite
dimensional to the infinite dimensional case.
\end{abstract}

\section{Preliminaries}

This paper is concerned with the study of approximate controllability of an
infinite dimensional stochastic equation with multiplicative noise 
\begin{equation}
\left\{ 
\begin{array}{l}
dX_{t}^{x,u}=\left( AX_{t}^{x,u}+Bu_{t}\right) dt+CX_{t}^{x,u}dW_{t}, \\ 
X_{0}=x\in H,%
\end{array}%
\right. 
\end{equation}%
where $u$ is a $U$-valued stochastic control process, and the state space $H$
as well as the control state space $U$ are separable real Hilbert spaces. We
say that the above equation enjoys the approximate controllability property
if, for any initial data $x\in H$, and all finite time horizon $T>0$, one
can find a control process $u$ which keeps the solution $X_{T}^{x,u}$
arbitrarily close to a given square integrable final condition.

For deterministic control systems with finite dimensional state space $%
\mathbb{C}
^{n}$, controllability is completely characterized by the well-known Kalman
condition. Often, it is convenient to study the observability of the adjoint
system rather than the controllability of the initial system. Indeed,
whenever dealing with a deterministic control system%
\begin{equation}
\left\{ 
\begin{array}{l}
dX_{t}^{x,u}=\left( AX_{t}^{x,u}+Bu_{t}\right) dt, \\ 
X_{0}=x\in 
\mathbb{C}
^{n},%
\end{array}%
\right.   \label{d1}
\end{equation}%
controllability is equivalent to the observability of the dual system%
\begin{equation}
\left\{ 
\begin{array}{l}
dY_{t}^{y}=-A^{\ast }Y_{t}^{y}dt,\text{ }O_{t}^{y}=B^{\ast }Y_{T}^{y}, \\ 
Y_{0}^{y}=y.%
\end{array}%
\right.   \label{d2}
\end{equation}%
A very powerful tool for this approach is the Hautus test. According to this
test, observability of (\ref{d2}) (and, thus, controllability for (\ref{d1}%
)) is equivalent to%
\begin{equation*}
rank\left[ 
\begin{array}{c}
sI-A^{\ast } \\ 
B^{\ast }%
\end{array}%
\right] =n,\text{ for all }s\in 
\mathbb{C}
.
\end{equation*}%
In the case of separable Hilbert state space, whenever $A$ generates an
exponentially stable semigroup, Russell and Weiss \cite{rw} have obtained a
necessary condition for observability which generalizes the Hautus
criterion. They have also conjectured that this condition is even
sufficient. Jacob and Zwart \cite{jz} proved that the above conjecture holds
true for the class of diagonal systems satisfying the strong stability
condition whenever the output space is finite dimensional. Similar arguments
allow to obtain in \cite{jp} a characterization of approximate
controllability of a deterministic controlled system with 1-dimensional
input.

In the stochastic framework, Kalman-type characterizations of \ approximate
controllability have been obtained, for the finite-dimensional case, by
Buckdahn, Quincampoix and Tessitore \cite{bqt} when the noise term is not
controlled, and by Goreac \cite{g} when the control is allowed to act on the
noise. The method they use relies on the duality between approximate
controllability and approximate observability for the dual equation. Riccati
algebraic arguments allow to obtain in \cite{bqt} and \cite{g} an invariance
criterion for the approximate controllability of the initial system.

In the case of controlled stochastic systems with infinite-dimensional state
space, we cite Barbu, R\u{a}\c{s}canu, Tessitore \cite{brt}, Fernandez-Cara,
Garrido-Atienza, Real \cite{fcgar}, and Sirbu, Tessitore \cite{st}. In \cite%
{st}, the authors characterize the property of (null) controllability with
the help of singular Riccati equations. They also provide a Riccati
characterization using the duality approach.

In this paper, we prove the duality between approximate controllability for
the forward system and some approximate observability for the dual system,
and we use this approach to show that the generalized Hautus test is a
necessary condition for approximate controllability whenever $A$ is the
generator of an exponentially stable semigroup.

The paper is organized as follows: In the first section we introduce the
standard notations and assumptions which will be used in what follows.
After, in the second section, we investigate the existence and the
uniqueness of the mild solution of the following backward stochastic
differential equation which is associated as dual equation to the controlled
system (1): 
\begin{equation*}
\left\{ 
\begin{array}{c}
dY_{t}=-\left( A^{\ast}Y_{t}+C^{\ast}Z_{t}\right) dt+Z_{t}dW_{t}, \\ 
Y_{T}=\xi\in L^{2}\left( \Omega,\mathcal{F}_{T},P;H\right) .%
\end{array}
\right.
\end{equation*}
We emphasize that the drift term in our dual backward equation contains not
only the unbounded operator $A^{\ast}$ acting on $Y$ but also the unbounded
operator $C^{\ast}$ that acts on $Z.$ To overcome the difficulties related
with, we make a joint dissipativity hypothesis which corresponds, in the
case of general heat equations, to the usual joint ellipticity condition.
Under these minimal assumptions we are able to prove the existence and the
uniqueness. Moreover, we provide a duality result between approximate
controllability for the forward equation and the approximate observability
of the dual system. The third section proves that, whenever $A$ generates an
exponentially stable semigroup, the Russell and Weiss generalization of the
Hautus test is a necessary condition for approximate controllability of
stochastic systems. Finally, we discuss \ as example the general heat
equation.

\section{Introduction}

Let us begin by introducing some basic notations and standard assumptions.
The spaces $\left( H,\left\langle \cdot,\cdot\right\rangle _{H}\right) ,$ $%
\left( U,\left\langle \cdot,\cdot\right\rangle _{U}\right) ,\left(
\Xi,\left\langle \cdot,\cdot\right\rangle _{\Xi}\right) $ are separable real
Hilbert spaces. We let $\mathcal{L}(\Xi,H)$ denote the space of all bounded $%
H$-valued linear operators on $\Xi,$ and $L_{2}(\Xi,H)$ be the subspace of
Hilbert-Schmidt operators. Both spaces are endowed with the usual norms.
Moreover, we consider a linear dissipative operator $A:D(A)\subset
H\longrightarrow H$ which generates a $C_{0}$-semigroup of linear operators $%
\left( e^{tA}\right) _{t\geq0},$ a linear bounded operator $B\in \mathcal{L}%
(U,H)$ and a linear operator $C:H\longrightarrow\mathcal{L}(\Xi,H)$ such
that, for all $t>0,$%
\begin{align*}
& a)\text{ }e^{tA}C\in\mathcal{L}\left( H;L_{2}(\Xi,H)\right) , \\
& b)\text{ }\left\vert e^{tA}C\right\vert _{\mathcal{L}\left(
H;L_{2}(\Xi,H)\right) }\leq Lt^{-\gamma},
\end{align*}
for some constants $\gamma\in\left[ 0,\frac{1}{2}\right) $ and $L>0.$

Let $\left( \Omega,\mathcal{F},P\right) $ be a complete probability space
endowed with a filtration $\left( \mathcal{F}_{t}\right) _{t\geq0}$ which is
supposed to satisfy the usual assumptions of completeness and
right-continuity. We denote by $W$ a cylindrical ($\mathcal{F}_{t}$)$-$%
Wiener process that takes its values in $\Xi.$ Finally, we let $\mathcal{U}$
denote the space of all $\left( \mathcal{F}_{t}\right) -$progressively
measurable processes $u:%
\mathbb{R}
_{+}\times\Omega\longrightarrow U$ such that%
\begin{equation*}
E\left[ \int_{0}^{T}\left\vert u_{t}\right\vert ^{2}dt\right] <\infty,\text{
for all }T>0.
\end{equation*}
$.$

The aim of this paper is to give an easy and verifiable criterion for
approximate controllability for the following linear stochastic differential
equation%
\begin{equation}
\left\{ 
\begin{array}{l}
dX_{t}^{x,u}=\left( AX_{t}^{x,u}+Bu_{t}\right) dt+CX_{t}^{x,u}dW_{t},\text{ }%
t\geq0. \\ 
X_{0}=x\in H.%
\end{array}
\right.  \label{eq1}
\end{equation}
Given an admissible control process $u\in\mathcal{U}$, an $\left( \mathcal{F}%
_{t}\right) $-progressively measurable process $X^{x,u}$ with 
\begin{equation*}
E\left[ \sup_{s\in\lbrack0,T]}\left\vert X_{s}^{x,u}\right\vert ^{2}\right]
<\infty,\text{ for all }T>0,
\end{equation*}
is a mild solution of (\ref{eq1}) if, for all $t>0,$%
\begin{equation}
X_{t}=e^{tA}x+\int_{0}^{t}e^{sA}Bu_{s}ds+\int_{0}^{t}e^{sA}CX_{s}dW_{s},
\label{eq2}
\end{equation}
$P$-a.s. Under the standard assumptions given above, there exists a unique
mild solution of (\ref{eq1}). For further results on mild solutions, the
reader is referred to Da Prato, Zabczyk \cite{dpz}, and Fuhrman, Tessitore 
\cite{ft}.

\section{The dual equation}

Let us now consider the following backward stochastic differential equation%
\begin{equation}
\left\{ 
\begin{array}{c}
dY_{t}=-\left( A^{\ast}Y_{t}+C^{\ast}Z_{t}\right) dt+Z_{t}dW_{t}, \\ 
Y_{T}=\xi\in L^{2}\left( \Omega,\mathcal{F}_{T},P;H\right) .%
\end{array}
\right.  \label{eq3}
\end{equation}
Since $C:H\longrightarrow\mathcal{L}(\Xi,H),$ also $Ce^{tA}:H\longrightarrow 
\mathcal{L}(\Xi,H)$, for all $t\geq0.$ Let us assume that, for all $t>0$,
all the values of $Ce^{tA}$ are in $L_{2}\left( \Xi;H\right) ,$ 
\begin{equation*}
Ce^{tA}:H\longrightarrow L_{2}\left( \Xi;H\right) .
\end{equation*}
Then, of course, the linear operator $\left( Ce^{tA}\right) ^{\ast}$ maps $%
L_{2}\left( \Xi;H\right) $ into $H$ and we can introduce the notion of a
mild solution for equation (\ref{eq3}). A mild solution of \ (\ref{eq3}) is
a couple $(Y,Z)$ of progressively measurable processes with values in $H,$
respectively $L_{2}(\Xi,H)$, such that 
\begin{equation*}
\left\{ 
\begin{array}{l}
(Y,Z)\in C\left( \left[ 0,T\right] ;L^{2}\left( \Omega;H\right) \right)
\times L^{2}\left( \left[ 0,T\right] \times\Omega;L_{2}\left( \Xi;H\right)
\right) , \\ 
\sup_{t\in\lbrack0,T]}E\left[ \left\vert Y_{t}\right\vert ^{2}\right] +E%
\left[ \int_{0}^{T}\left\vert Z_{t}\right\vert ^{2}dt\right] <\infty, \\ 
\int_{0}^{T}\left\vert \left( Ce^{\left( s-t\right) A}\right)
^{\ast}Z_{s}\right\vert ds<\infty,\text{ }P-a.s., \\ 
\\ 
Y_{t}=e^{(T-t)A^{\ast}}\xi+\int_{t}^{T}\left( Ce^{(s-t)A}\right) ^{\ast
}Z_{s}ds-\int_{t}^{T}e^{(s-t)A^{\ast}}Z_{s}dW_{s},\text{ }t\in\left[ 0,T%
\right] .%
\end{array}
\right.
\end{equation*}
If $C$ is a bounded linear operator, then it has been shown in Confortola [%
\cite{c}, Th. 2.2] that (\ref{eq3}) admits a unique mild solution. Let us
suppose that

\textbf{(A1)}\textit{\ The operator }$C$\textit{\ may be written as sum of
two linear operators }$C_{1},$\textit{\ }$C_{2}$\textit{\ }%
\begin{equation*}
C=C_{1}+C_{2},
\end{equation*}
\textit{satisfying the following properties:}

\textit{1) }$C_{2}$\textit{\ is a bounded operator from }$H$\textit{\ to }$%
L_{2}\left( \Xi;H\right) $\textit{,}

\textit{2) for all }$t>0,$\textit{\ }$C_{1}e^{tA}\in L\left( H;L_{2}\left(
\Xi;H\right) \right) .$\textit{\ Moreover, we suppose that there exist some }%
$\gamma\in\left[ 0,\frac{1}{2}\right) $\textit{\ and some positive constant }%
$L>0$\textit{\ such that }%
\begin{equation*}
\left\vert C_{1}e^{tA}\right\vert _{\mathcal{L}\left( H;L_{2}\left(
\Xi;H\right) \right) }\leq Lt^{-\gamma},
\end{equation*}
\textit{for all }$t>0.$

\textit{3) There exists some constant }$a>\frac{1}{2}$\textit{\ such that }%
\begin{equation*}
\text{ }A+a\left( C_{1}e^{\delta A}\right) ^{\ast}\left( C_{1}e^{\delta
A}\right) \text{ is dissipative},
\end{equation*}
\textit{for some sequence }$\delta\searrow0.$

If $C_{1}$ is different of zero, we shall also assume that

\textbf{(A2)} $-A^{2}$\textit{\ is dissipative.}

\begin{remark}
\textit{If }$A$\textit{\ is a self-adjoint, dissipative operator which
generates a contraction semigroup, then (A.2) is obviously satisfied.}

\textit{Moreover, if we suppose that }$C_{1}$\textit{\ takes its values in }$%
L_{2}\left( \Xi;H\right) ,$\textit{\ then we may replace (A1) 3) by}

\textit{3') there exists some constant }$a>\frac{1}{2}$\textit{\ such that }%
\begin{equation*}
\text{ }A+aC_{1}^{\ast}C_{1}\text{ is dissipative.}
\end{equation*}

\textit{Indeed, in this case }$e^{2\delta A}$\textit{\ is a bounded operator
which commutes with the self-adjoint positive operator }$-A$\textit{\ and
also with its square root }$\sqrt{-A}.$ \textit{\ Thus, for all }$x\in D(A),$%
\begin{equation*}
\left\langle e^{2\delta A}(-A)x,x\right\rangle =\left\vert e^{\delta A}\sqrt{%
-A}x\right\vert ^{2}\leq\left\vert \sqrt{-A}x\right\vert ^{2}=\left\langle
(-A)x,x\right\rangle .
\end{equation*}
\textit{It follows that }$A-e^{\delta A}Ae^{\delta A^{\ast}}$\textit{\ is
dissipative. Therefore, also }$A+ae^{\delta
A^{\ast}}C_{1}^{\ast}C_{1}e^{\delta A}$\textit{\ is dissipative.}
\end{remark}

We now can state the main result of this section.

\begin{theorem}
Under the assumptions (A1) and (A2), there exists a unique mild solution of
the backward linear stochastic differential equation (\ref{eq3}). Moreover,
this solution satisfies 
\begin{equation}
\sup_{t\in\lbrack0,T]}E\left[ \left\vert Y_{t}\right\vert ^{2}\right] +E%
\left[ \int_{0}^{T}\left\vert Z_{s}\right\vert ^{2}ds\right] \leq kE\left[
\left\vert \xi\right\vert ^{2}\right] ,  \label{ineq0.4}
\end{equation}
where $k>0$ is some constant that doesn't depend on the particular choice of 
$\xi$ but only on the operators $A,C$ and the time horizon $T$.
\end{theorem}

\begin{remark}
\textit{1. The existence and uniqueness of the solution for equation (\ref%
{eq3}) has been studied by Tessitore \cite{t} for the case in which }$A$%
\textit{\ generates an analytic semigroup of contractions of negative type;
the Brownian motion was supposed to be finite-dimensional. His main
assumption, the joint dissipativity condition, was justified by its
necessity for the "well-posedness" and coercivity of the forward system. The
approach is fundamentally different from ours and relies on duality methods.
However, let us point out that the author obtains, for his analytic case,
stronger space regularity properties for the solution of the BSDE.}

\textit{2. Ma, Yong \cite{my2} treated a particular linear, degenerate
BSPDE. Their method relies on a parabolicity assumption and a priori
estimates that allowed the authors to get the well-posedness of the problem,
the existence, the uniqueness as well as regularity properties. Later the
same technique was used by Hu, Ma, Yong \cite{hmy} for further extensions.}
\end{remark}

\begin{proof}
(of Theorem 1). We begin by proving the \underline{existence}: The main
difficulty to prove the existence and the uniqueness for a BSDE in infinite
dimensions with unbounded linear operators consists in the fact that It\^{o}%
's formula can't be applied directly to this equation because it is defined
only in the mild sense. To overcome this difficulty, we have to reduce the
problem with the help of two different approximations to BSDEs that allow
the application of It\^{o}'s formula. We first approximate our original BSDE
by the following one: 
\begin{equation}
\left\{ 
\begin{array}{l}
dY_{t}^{\delta}=-A^{\ast}Y_{t}^{\delta}dt-(C_{1}e^{\delta
A})^{\ast}Z_{t}^{\delta}dt-C_{2}^{\ast}Z_{t}^{\delta}dt+Z_{t}^{\delta}dW_{t},
\\ 
Y_{T}^{\delta}=\xi\in L^{2}(\Omega,\mathcal{F}_{T},P;H)%
\end{array}
\right.  \label{eqdelta}
\end{equation}
For this approximating equation we know that, due to the results of
Confortola \cite{c}, there exists a unique mild solution $%
(Y^{\delta},Z^{\delta})$ for every $\delta>0$.

\smallskip

In a first step we prove that

\medskip

\underline{Step 1}. There is a positive constant $k$ independent of $%
\delta>0 $ and $\xi$ such that

\begin{equation}
\sup_{t\in\lbrack0,T]}E\left[ \left\vert Y_{t}^{\delta}\right\vert ^{2}%
\right] +E\left[ \int_{0}^{T}\left\vert Z_{s}^{\delta}\right\vert ^{2}ds%
\right] \leq kE\left\vert \xi\right\vert ^{2}.  \label{ineg2}
\end{equation}

\medskip

Indeed, we introduce the Yosida approximation of the dissipative operator $%
A^{\ast }$, $A_{n}^{\ast }=n(nI-A^{\ast })^{-1}A^{\ast }=J_{n}^{\ast
}A^{\ast }$, and we consider the following approximating BSDE:%
\begin{equation*}
\left\{ 
\begin{array}{l}
dY_{t}^{n,\delta }=-\ A_{n}^{\ast }Y_{t}^{n,\delta }dt-J_{n}^{\ast }\left(
C_{1}e^{\delta A}\right) ^{\ast }Z_{t}^{n,\delta }dt-C_{2}^{\ast
}Z_{t}^{n,\delta }dt+Z_{t}^{n,\delta }dW_{t}, \\ 
Y_{T}^{n,\delta }=\xi \in L^{2}(\Omega ,\mathcal{F}_{T},P;H).%
\end{array}%
\right.
\end{equation*}%
It is well known that the above equation admits a unique solution $\left(
Y^{n,\delta },Z^{n,\delta }\right) $. Let $1<\alpha <2a$ and $\beta >0$ be
such that $\frac{1}{\alpha }+\frac{1}{\beta }<1.$ Then, by applying It\^{o}%
's formula to $\left\vert Y^{n,\delta }\right\vert ^{2}$ we obtain 
\begin{align}
E\left\vert \xi \right\vert ^{2}& =E\left[ \left\vert Y_{t}^{n,\delta
}\right\vert ^{2}\right] -2E\left[ \int_{t}^{T}\left\langle A_{n}^{\ast
}Y_{s}^{n,\delta },Y_{s}^{n,\delta }\right\rangle \right]  \notag \\
& -2E\left[ \int_{t}^{T}\left\langle J_{n}^{\ast }\left( C_{1}e^{\delta
A}\right) ^{\ast }Z_{s}^{n,\delta },Y_{s}^{n,\delta }\right\rangle \right] 
\notag \\
& -2E\left[ \int_{t}^{T}\left\langle C_{2}^{\ast }Z_{s}^{n,\delta
},Y_{s}^{n,\delta }\right\rangle \right] +E\left[ \int_{t}^{T}\left\vert
Z_{s}^{n,\delta }\right\vert ^{2}ds\right]  \notag \\
& \geq E\left[ \left\vert Y_{t}^{n,\delta }\right\vert ^{2}\right] +\left( 1-%
\frac{1}{\alpha }-\frac{1}{\beta }\right) E\left[ \int_{t}^{T}\left\vert
Z_{s}^{n,\delta }\right\vert ^{2}ds\right]  \notag \\
& -2E\left[ \int_{t}^{T}\left\langle \left( A_{n}^{\ast }+\frac{\alpha }{2}%
J_{n}^{\ast }\left( C_{1}e^{\delta A}\right) ^{\ast }\left( C_{1}e^{\delta
A}\right) J_{n}\right) Y_{s}^{n,\delta },Y_{s}^{n,\delta }\right\rangle %
\right]  \notag \\
& -\beta |C_{2}^{\ast }|^{2}E\left[ \int_{t}^{T}\left\vert Y_{s}^{n,\delta
}\right\vert ^{2}ds\right] ,  \label{eg1}
\end{align}%
On the other hand, with the help of assumption (A.2) we can prove that 
\begin{align*}
& A_{n}^{\ast }+\frac{\alpha }{2}J_{n}^{\ast }\left( C_{1}e^{\delta
A}\right) ^{\ast }\left( C_{1}e^{\delta A}\right) J_{n} \\
& =-n^{-1}A_{n}^{\ast }A_{n}+J_{n}^{\ast }\left( A^{\ast }+\frac{\alpha }{2}%
\left( C_{1}e^{\delta A}\right) ^{\ast }\left( C_{1}e^{\delta A}\right)
\right) J_{n}
\end{align*}%
is a dissipative operator. It then follows from (\ref{eg1}) that 
\begin{align*}
& E\left[ \left\vert Y_{t}^{n,\delta }\right\vert ^{2}\right] +\left( 1-%
\frac{1}{\alpha }-\frac{1}{\beta }\right) E\left[ \int_{t}^{T}\left\vert
Z_{s}^{n,\delta }\right\vert ^{2}ds\right] \\
& \leq E\left\vert \xi \right\vert ^{2}+\beta |C_{2}|^{2}E\left[
\int_{t}^{T}\left\vert Y_{s}^{n,\delta }\right\vert ^{2}ds\right] ,
\end{align*}%
and Gronwall's inequality yields%
\begin{equation}
\sup_{t\in \lbrack 0,T]}E\left[ \left\vert Y_{t}^{n,\delta }\right\vert ^{2}%
\right] +E\left[ \int_{0}^{T}\left\vert Z_{s}^{n,\delta }\right\vert ^{2}ds%
\right] \leq kE\left\vert \xi \right\vert ^{2},  \label{ineg1}
\end{equation}%
Notice that the constant $k$ here is independent of $n\geq 1,\delta >0$ and
of $\xi $; it denotes a generic constant whose value can change from line to
line. From the above estimate we can conclude that there is a subsequence,
still denoted $\left( Y^{n,\delta },Z^{n,\delta }\right) _{n}$, such that $%
Y^{n,\delta }\rightarrow Y^{\delta }$ weakly * in $L^{\infty }\left(
[0,T];L^{2}(\Omega ;H)\right) $ and $Z^{n,\delta }\rightarrow Z^{\delta }$
weakly in $L^{2}\left( \Omega \times \lbrack 0,T];L_{2}(\Xi ;H)\right) .$ It
can be easily proved the limit $(Y^{\delta },Z^{\delta })$ is the unique
mild solution of\ (\ref{eqdelta}). This allows to consider for $Y^{\delta }$
its version in $C\left( [0,T];L^{2}(\Omega ;H)\right) .$ Finally, from
Mazur's theorem we obtain that $(Y^{\delta },Z^{\delta })$ satisfies the
estimate announced in step 1.

\medskip

In preparation of the next step we observe that, since $(Y^{\delta
},Z^{\delta })_{\delta >0}$ is bounded in $L^{\infty }\left(
[0,T];L^{2}(\Omega ;H)\right) \times L^{2}\left( \Omega \times \lbrack
0,T];L_{2}(\Xi ;H)\right) ,$ we get the existence of some subsequence, again
denoted by $(Y^{\delta },Z^{\delta })_{\delta >0}$, such that $Y^{\delta
}\rightarrow Y$ weak * in $L^{\infty }\left( [0,T];L^{2}(\Omega ;H)\right) $
and $Z^{\delta }\rightarrow Z$ weakly in \newline
$L^{2}\left( \Omega \times \lbrack 0,T];L_{2}(\Xi ;H)\right) $, as $\delta
\rightarrow 0.$

\medskip

We want to prove that the couple $\left( Y,Z\right) $ obtained above is a
mild solution of our BSDE:\ 

\begin{align}
Y_{t} & =e^{(T-t)A^{\ast}}\xi+\int_{t}^{T}\left( C_{1}e^{\left( s-t\right)
A}\right) ^{\ast}Z_{s}ds  \notag \\
&
+\int_{t}^{T}e^{(s-t)A^{\ast}}C_{2}^{\ast}Z_{s}ds-\int_{t}^{T}e^{(s-t)A^{%
\ast}}Z_{s}dW_{s}.  \label{eqmild}
\end{align}

For this we notice that, since $(Y^{\delta},Z^{\delta})$ is a mild solution
of (\ref{eqdelta}), we have%
\begin{align}
Y_{t}^{\delta} & =e^{(T-t)A^{\ast}}\xi+\int_{t}^{T}e^{(s-t)A^{\ast}}\left(
C_{1}e^{\delta A}\right) ^{\ast}Z_{s}^{\delta}ds  \notag \\
&
+\int_{t}^{T}e^{(s-t)A^{\ast}}C_{2}^{\ast}Z_{s}^{\delta}ds-%
\int_{t}^{T}e^{(s-t)A^{\ast}}Z_{s}^{\delta}dW_{s}  \label{eqdeltamild}
\end{align}
and we show the following:

\medskip

\underline{Step 2} The process 
\begin{equation*}
M_{t}^{1,\delta }=\int_{t}^{T}e^{(s-t)A^{\ast }}\left( C_{1}e^{\delta
A}\right) ^{\ast }Z_{s}^{\delta }ds,\,t\in \lbrack 0,T],
\end{equation*}%
belongs to $L^{\infty }\left( [0,T];L^{2}(\Omega ;H)\right) $ and converges
weakly * in \newline
$L^{\infty }\left( [0,T];L^{2}(\Omega ;H)\right) $ to $M^{1}=\left(
\int_{t}^{T}\left( C_{1}e^{\left( s-t\right) A}\right) ^{\ast
}Z_{s}ds\right) _{t\in \lbrack 0,T]}$.

\medskip

Indeed, by using that 
\begin{equation*}
e^{\delta ^{\prime }A^{\ast }}\left( C_{1}e^{\delta A}\right) ^{\ast
}=\left( C_{1}e^{\left( \delta +\delta ^{\prime }\right) A}\right) ^{\ast },
\end{equation*}%
for all $\delta ,\delta ^{^{\prime }}>0$, we have 
\begin{align*}
& E\left[ \left\vert \int_{t}^{T}e^{(s-t)A^{\ast }}\left( C_{1}e^{\delta
A}\right) ^{\ast }Z_{s}^{\delta }ds\right\vert ^{2}\right] \\
& \leq E\left[ \left( \int_{t}^{T}\left\vert e^{\delta A^{\ast }}\left(
C_{1}e^{(s-t)A}\right) ^{\ast }Z_{s}^{\delta }\right\vert ds\right) ^{2}%
\right] \\
& \leq kE\left[ \int_{t}^{T}\left( s-t\right) ^{-2\gamma
}ds\int_{t}^{T}\left\vert Z_{s}^{\delta }\right\vert ^{2}ds\right] \\
& \leq kE\left\vert \xi \right\vert ^{2},
\end{align*}%
which implies that $\left\{ M^{1,\delta },\delta >0\right\} \subset
L^{\infty }\left( [0,T];L^{2}(\Omega ;H)\right) $ is bounded. Moreover, for
all $\phi \in L^{2}\left( \Omega ;H\right) $ and $t\in \lbrack 0,T]$, 
\begin{align}
E\left[ \left\langle M_{t}^{1,\delta },\phi \right\rangle \right] & =E\left[
\int_{t}^{T}\left\langle \left( C_{1}e^{(s-t)A}\right) ^{\ast }Z_{s}^{\delta
},\left( e^{\delta A^{\ast }}-I\right) \phi \right\rangle ds\right]  \notag
\\
& +E\left[ \int_{t}^{T}\left\langle \left( C_{1}e^{(s-t)A}\right) ^{\ast
}Z_{s}^{\delta },\phi \right\rangle ds\right] =:I_{1}^{\delta
}+I_{2}^{\delta },  \label{eg2}
\end{align}%
where 
\begin{align*}
I_{1}^{\delta }& =E\left[ \left\vert \int_{t}^{T}\left\langle \left(
C_{1}e^{(s-t)A}\right) ^{\ast }Z_{s}^{\delta },\left( e^{\delta A^{\ast
}}-I\right) \phi \right\rangle ds\right\vert \right] \\
& \leq E\left[ \left\vert \left( e^{\delta A^{\ast }}-I\right) \phi
\right\vert \int_{t}^{T}\left\vert \left( C_{1}e^{(s-t)A}\right) ^{\ast
}Z_{s}^{\delta }\right\vert ds\right] \\
& \leq \left( E\left[ \int_{t}^{T}\left( s-t\right) ^{-2\gamma
}ds\int_{t}^{T}\left\vert Z_{s}^{\delta }\right\vert ^{2}ds\right] \right) ^{%
\frac{1}{2}}\left( E\left[ \left\vert \left( e^{\delta A^{\ast }}-I\right)
\phi \right\vert ^{2}\right] \right) ^{\frac{1}{2}} \\
& \leq k\left( E\left\vert \xi ^{2}\right\vert \right) ^{\frac{1}{2}}\left( E%
\left[ \left\vert \left( e^{\delta A^{\ast }}-I\right) \phi \right\vert ^{2}%
\right] \right) ^{\frac{1}{2}}.
\end{align*}%
Consequently, due to the dominated convergence theorem, 
\begin{equation*}
I_{1}^{\delta }\rightarrow 0\text{ as }\delta \rightarrow 0.
\end{equation*}%
For the second term we have 
\begin{equation*}
I_{2}^{\delta }=E\left[ \int_{t}^{T}\left\langle \left(
C_{1}e^{(s-t)A}\right) ^{\ast }Z_{s}^{\delta },\phi \right\rangle ds\right]
=E\left[ \int_{t}^{T}\left\langle Z_{s}^{\delta },\left(
C_{1}e^{(s-t)A}\right) \phi \right\rangle ds\right] ,
\end{equation*}%
and since $\left( C_{1}e^{(s-t)A}\right) \phi \in L^{2}\left( \Omega \times
\lbrack 0,T];L_{2}(\Xi ;H)\right) ,$ it follows from the weak convergence of 
$Z^{\delta }$ to $Z$ that%
\begin{equation*}
I_{2}^{\delta }=E\left[ \int_{t}^{T}\left\langle \left(
C_{1}e^{(s-t)A}\right) ^{\ast }Z_{s}^{\delta },\phi \right\rangle ds\right]
\rightarrow E\left[ \int_{t}^{T}\left\langle \left( C_{1}e^{(s-t)A}\right)
^{\ast }Z_{s},\phi \right\rangle ds\right] ,
\end{equation*}%
and from (\ref{eg2}) we then get%
\begin{equation*}
E\left[ \left\langle M_{t}^{1,\delta },\phi \right\rangle \right]
\rightarrow E\left[ \left\langle M_{t}^{1},\phi \right\rangle \right] \text{
as }\delta \rightarrow 0.
\end{equation*}%
In order to prove that $M^{1,\delta }$ converges in the weak * topology on $%
L^{\infty }([0,T];$ $L^{2}(\Omega ;H))$ to $M^{1},$ we consider $\Phi \in
L^{1}\left( [0,T];L^{2}(\Omega ;H)\right) $, and use the fact that, for all $%
t\in \lbrack 0,T]$ for which $\Phi _{t}\in L^{2}(\Omega ;H),$ the previous
convergence holds with $\Phi _{t}$ at the place of $\phi $. We then apply a
dominated convergence argument and get the statement of step 2.

\medskip

\underline{Step 3}. The couple $(Y,Z)$ is a solution of the BSDE 
\begin{align}
Y_{t} & =e^{(T-t)A^{\ast}}\xi+\int_{t}^{T}\left( C_{1}e^{(s-t)A}\right)
^{\ast}Z_{s}ds  \notag \\
&
+\int_{t}^{T}e^{(s-t)A^{\ast}}C_{2}^{\ast}Z_{s}ds-\int_{t}^{T}e^{(s-t)A^{%
\ast}}Z_{s}dW_{s}.  \label{star}
\end{align}
Moreover, 
\begin{equation}
\sup_{t\in\lbrack0,T]}E\left[ \left\vert Y_{t}\right\vert ^{2}\right] +E%
\left[ \int_{0}^{T}\left\vert Z_{s}\right\vert ^{2}ds\right] \leq
kE\left\vert \xi\right\vert ^{2}.
\end{equation}

\medskip

To prove the above statement we write $Y_{t}^{\delta },\,t\in \lbrack 0,T],$
as 
\begin{equation*}
Y_{t}^{\delta }=e^{(T-t)A^{\ast }}\xi +M_{t}^{1,\delta }+M_{t}^{2,\delta
}+M_{t}^{3,\delta }.
\end{equation*}%
While we have already studied the convergence of $M^{1,\delta }$ in the
preceding step, it is an immediate consequence of the boundedness of the
operator $C_{2}$ that $M_{t}^{2,\delta }=\int_{t}^{T}e^{(s-t)A^{\ast
}}C_{2}^{\ast }Z_{s}^{\delta }ds$ converges weakly * in $L^{\infty }\left(
[0,T];L^{2}(\Omega ;H)\right) $ to $M_{t}^{2}=\int_{t}^{T}e^{(s-t)A^{\ast
}}C_{2}^{\ast }Z_{s}ds.$

For the noise term $M_{t}^{3,\delta }=\int_{t}^{T}e^{(s-t)A^{\ast
}}Z_{s}^{\delta }dW_{s}$ we notice that since $Z^{\delta }$ converges weakly
in $L^{2}\left( \Omega \times \lbrack 0,T];L_{2}(\Xi ;H)\right) $ to $Z,$ $%
e^{(\cdot -t)A^{\ast }}Z_{\cdot }^{\delta }$ also converges weakly to $%
e^{(\cdot -t)A^{\ast }}Z_{\cdot }.$ We apply the martingale representation
theorem to get that \newline
$\int_{t}^{T}e^{(s-t)A^{\ast }}Z_{s}^{\delta }dW_{s}$ converges weakly in $%
L^{2}\left( \Omega ;H\right) $ to $\int_{t}^{T}e^{(s-t)A^{\ast }}Z_{s}dW_{s}.
$ Using, as before, the dominated convergence, we get that%
\begin{eqnarray*}
&&N_{t}^{\delta }=\int_{t}^{T}e^{(s-t)A^{\ast }}Z_{s}^{\delta }dW_{s}\text{
converges in the weak* topology on } \\
&&L^{\infty }\left( [0,T];L^{2}(\Omega ;H)\right) \text{ to }%
N_{t}=\int_{t}^{T}e^{(s-t)A^{\ast }}Z_{s}dW_{s}.
\end{eqnarray*}

We now pass to the $L^{\infty}\left( [0,T];L^{2}(\Omega;H)\right) $ weak *
limit in the approximating mild equation (\ref{eqdeltamild}). This yields
the statement of step 3, with the only difference, that for the BSDE which
has been got by a weak limit, we only know for the moment that this equation
is satisfied $dtdP$-a.e. To obtain that the BSDE is satisfied by $(Y,Z)$ for
all time points of the interval $[0,T]$, $P$-a.s., we need the following
auxiliary statement:
\end{proof}

\begin{lemma}
The process

$\Phi_{t}=e^{\left( T-t\right) A^{\ast}}\xi+\left( C_{1}e^{(r-t)A}\right)
^{\ast}Z_{r}dr+\int_{t}^{T}e^{\left( r-t\right) A^{\ast}}C_{2}^{\ast}Z_{r}dr$

$\qquad\qquad-\int_{t}^{T}e^{\left( r-t\right) A^{\ast}}Z_{r}dW_{r},\, t\in[%
0,T]$, is mean-square continuous.
\end{lemma}

\begin{proof}
We return to the proof of our theorem. The proof of the lemma will be given
afterwards.

The above result allows to conclude the proof of step 3. Indeed, the above
lemma guarantees the existence of a version of the solution $\left(
Y,Z\right) $ in $C\left( \left[ 0,T\right] ;L^{2}\left( \Omega;H\right)
\right) \times$ $L^{2}\left( \Omega\times\lbrack0,T];L_{2}(\Xi;H)\right) .$
For this version we have (\ref{star}) for all $t\in\left[ 0,T\right] $.

\medskip

Let us prove now the \underline{uniqueness} of the solution of our BSDE. In
virtue of the linearity of the equation it suffices to prove the following:

\medskip

\underline{Step 4}. The only solution $(Y,Z)$ of the BSDE 
\begin{equation*}
\left\{ 
\begin{array}{l}
dY_{t}=-A^{\ast}Y_{t}dt-C^{\ast}Z_{t}dt+Z_{t}dW_{t}, \\ 
Y_{T}=0.%
\end{array}
\right.
\end{equation*}
is the trivial one: $(Y,Z)=(0,0)$.

To prove this, we have to transform the BSDE into an equation which allows
to apply It\^{o}'s formula. For this reason we put, for all $n\geq 1$ and $%
\delta >0$,%
\begin{equation*}
\widetilde{Y}_{\cdot }:=J_{n}^{\ast }e^{\delta A^{\ast }}Y_{\cdot },
\end{equation*}%
and we observe that the such introduced process $\widetilde{Y}$ satisfies
the following backward equation:%
\begin{equation*}
\left\{ 
\begin{array}{l}
d\widetilde{Y}_{t}=-A^{\ast }\widetilde{Y}_{t}dt-J_{n}^{\ast }\left(
C_{1}e^{\delta A}\right) ^{\ast }Z_{t}dt-J_{n}^{\ast }e^{\delta A^{\ast
}}C_{2}^{\ast }Z_{t}dt+J_{n}^{\ast }e^{\delta A^{\ast }}Z_{t}dW_{t}, \\ 
\widetilde{Y}_{T}=0.%
\end{array}%
\right.
\end{equation*}%
To this equation we can apply It\^{o}'s formula (Indeed, notice that $%
A^{\ast }\widetilde{Y}_{\cdot }=(J_{n}^{\ast }e^{\delta A^{\ast }}A^{\ast
})Y_{\cdot }$, where the operator $J_{n}^{\ast }e^{\delta A^{\ast }}A^{\ast
} $ is bounded). This yields:%
\begin{align}
& 0=E\left[ \left\vert J_{n}^{\ast }e^{\delta A^{\ast }}Y_{t}\right\vert ^{2}%
\right] -2E\left[ \int_{t}^{T}\left\langle A^{\ast }\widetilde{Y}_{s},%
\widetilde{Y}_{s}\right\rangle ds\right]  \notag \\
& -2E\left[ \int_{t}^{T}\left\langle J_{n}^{\ast }\left( C_{1}e^{\delta
A}\right) ^{\ast }Z_{s},\widetilde{Y}_{s}\right\rangle ds\right]  \notag \\
& -2E\left[ \int_{t}^{T}\left\langle J_{n}^{\ast }e^{\delta A^{\ast
}}C_{2}^{\ast }Z_{s},\widetilde{Y}_{s}\right\rangle ds\right] +E\left[
\int_{t}^{T}\left\vert J_{n}^{\ast }e^{\delta A^{\ast }}Z_{s}\right\vert
^{2}ds\right]  \notag \\
& \geq E\left[ \left\vert J_{n}^{\ast }e^{\delta A^{\ast }}Y_{t}\right\vert
^{2}\right] -2E\left[ \int_{t}^{T}\left\langle \left( A^{\ast }+\frac{\alpha 
}{2}J_{n}^{\ast }\left( C_{1}e^{\delta A}\right) ^{\ast }\left(
C_{1}e^{\delta A}\right) J_{n}\right) \widetilde{Y}_{s},\widetilde{Y}%
_{s}\right\rangle ds\right]  \notag \\
& -\beta \left\vert C_{2}\right\vert ^{2}E\left[ \int_{t}^{T}\left\vert
Y_{s}\right\vert ^{2}ds\right] +E\left[ \int_{t}^{T}\left\vert J_{n}^{\ast
}e^{\delta A^{\ast }}Z_{s}\right\vert ^{2}ds\right]  \notag \\
& -\left( \frac{1}{\alpha }+\frac{1}{\beta }\right) E\left[
\int_{t}^{T}\left\vert Z_{s}\right\vert ^{2}ds\right] ,  \label{ineq0.6}
\end{align}%
To be able to go ahead with the above estimate we need the dissipativity of
the operator $A^{\ast }+\frac{\alpha }{2}J_{n}^{\ast }\left( C_{1}e^{\delta
A}\right) ^{\ast }\left( C_{1}e^{\delta A}\right) J_{n}$.

For this end we notice that 
\begin{equation*}
\left( nI-A^{\ast }\right) A^{\ast }\left( nI-A\right) -n^{2}A^{\ast
}=-nA^{\ast }A^{\ast }-nA^{\ast }A+A^{\ast }A^{\ast }A
\end{equation*}%
and apply this relation to the operator $\left( nI-A\right) ^{-1}.$ To the
relation we then apply $\left( nI-A^{\ast }\right) ^{-1}$. So we obtain the
following equality: 
\begin{equation*}
A^{\ast }-J_{n}^{\ast }A^{\ast }J_{n}=-n^{-1}J_{n}^{\ast }\left( A^{\ast
}\right) ^{2}J_{n}-n^{-1}J_{n}^{\ast }A^{\ast }AJ_{n}+n^{-2}J_{n}^{\ast
}A^{\ast }A^{\ast }AJ_{n},
\end{equation*}%
which proves that the operator $A^{\ast }-J_{n}^{\ast }A^{\ast }J_{n}$ is
dissipative. It now follows easily that also the operator%
\begin{align*}
& A^{\ast }+\frac{\alpha }{2}J_{n}^{\ast }\left( C_{1}e^{\delta A}\right)
^{\ast }\left( C_{1}e^{\delta A}\right) J_{n} \\
& =A^{\ast }-J_{n}^{\ast }A^{\ast }J_{n}+J_{n}^{\ast }A^{\ast }J_{n}+\frac{%
\alpha }{2}J_{n}^{\ast }\left( C_{1}e^{\delta A}\right) ^{\ast }\left(
C_{1}e^{\delta A}\right) J_{n}
\end{align*}%
is dissipative if the parameters $\alpha ,\beta $ are chosen as in (\ref{eg1}%
).

This dissipativity allows to go ahead in (\ref{ineq0.6}) and to conclude that%
\begin{align*}
& E\left[ \left\vert J_{n}^{\ast}e^{\delta A^{\ast}}Y_{t}\right\vert ^{2}%
\right] +E\left[ \int_{t}^{T}\left\vert J_{n}^{\ast}e^{\delta
A^{\ast}}Z_{s}\right\vert ^{2}ds\right] \\
& \leq\beta\left\vert C_{2}\right\vert ^{2}E\left[ \int_{t}^{T}\left\vert
Y_{s}\right\vert ^{2}ds\right] +\left( \frac{1}{\alpha}+\frac{1}{\beta }%
\right) E\left[ \int_{t}^{T}\left\vert Z_{s}\right\vert ^{2}ds\right] .
\end{align*}
Recall that $\left( Y,Z\right) \in$ $L^{2}\left( \Omega\times\left[ 0,T%
\right] ;H\times L_{2}\left( \Xi;H\right) \right) $. Thus, letting $%
n\rightarrow\infty$ and then $\delta\rightarrow0$ in the above estimate, we
get 
\begin{equation*}
E\left[ \left\vert Y_{t}\right\vert ^{2}\right] +\left( 1-\frac{1}{\alpha }-%
\frac{1}{\beta}\right) E\left[ \int_{t}^{T}\left\vert Z_{s}\right\vert ^{2}ds%
\right] \leq k\beta\left\vert C_{2}\right\vert ^{2}E\left[ \int
_{t}^{T}\left\vert Y_{s}\right\vert ^{2}ds\right] .
\end{equation*}
Finally, we take the supremum over $t\in\lbrack0,T]$ and apply Gronwall's
inequality. Thus we obtain 
\begin{equation*}
\sup_{t\in\lbrack0,T]}E\left[ \left\vert Y_{t}\right\vert ^{2}\right] +E%
\left[ \int_{0}^{T}\left\vert Z_{s}\right\vert ^{2}ds\right] =0,
\end{equation*}
and the claimed uniqueness follows as immediate consequence.$_{\blacksquare}$
\end{proof}

\medskip

In order to really complete the proof of the theorem we still have to give
the proof of Lemma 1.

\begin{proof}
(of Lemma 1) A standard estimate for the process $\Phi$ defined in Lemma 1
gives the following for all $s,t\geq0$:%
\begin{align}
E\left[ \left\vert \Phi_{t}-\Phi_{s}\right\vert ^{2}\right] & \leq k\left( E%
\left[ \left\vert \left( e^{\left\vert t-s\right\vert A^{\ast}}-I\right)
\Phi_{t\vee s}\right\vert ^{2}\right] \right.  \notag \\
& +E\left[ \left\vert \int_{s\wedge t}^{s\vee t}\left( C_{1}e^{\left(
r-s\right) A}\right) ^{\ast}Z_{r}\right\vert ^{2}\right]  \notag \\
& \left. +E\left[ \left\vert \int_{s\wedge t}^{s\vee t}e^{\left( r-s\right)
A^{\ast}}C_{2}^{\ast}Z_{r}\right\vert ^{2}\right] +E\left[ \int_{s\wedge
t}^{s\vee t}\left\vert e^{\left( r-t\right) A^{\ast}}Z_{r}\right\vert ^{2}dr%
\right] \right)  \notag \\
& \leq k\left( E\left[ \left\vert \left( e^{\left\vert t-s\right\vert
A^{\ast}}-I\right) \Phi_{t\vee s}\right\vert ^{2}\right] \right.  \notag \\
& \left. +\left( 1+\left\vert t-s\right\vert ^{1-2\gamma}\right) E\left[
\int_{s\wedge t}^{s\vee t}\left\vert Z_{r}\right\vert ^{2}dr\right] \right) .
\label{maj}
\end{align}
Here $k$ denotes a generic constant that is independent of $%
s,t\in\lbrack0,T] $ and can change from line to line.

Since $Z\in L^{2}\left( \Omega\times\left[ 0,T\right] ;L_{2}\left(
\Xi;H\right) \right) ,$ it is a direct consequence of the dominated
convergence theorem that%
\begin{equation*}
\lim_{s\rightarrow t}E\left[ \int_{s\wedge t}^{s\vee t}\left\vert
Z_{r}\right\vert ^{2}dr\right] =0.
\end{equation*}
It remains to show that also $E\left[ \left\vert \left( e^{\left\vert
t-s\right\vert A^{\ast}}-I\right) \Phi_{t\vee s}\right\vert ^{2}\right] $
converges to zero, as $s\rightarrow t$. We first consider this limit for $%
t>s\uparrow t$. In this case 
\begin{equation*}
E\left[ \left\vert \left( e^{\left\vert t-s\right\vert A^{\ast}}-I\right)
\Phi_{t\vee s}\right\vert ^{2}\right] =E\left[ \left\vert \left( e^{\left(
t-s\right) A^{\ast}}-I\right) \Phi_{t}\right\vert ^{2}\right] ,
\end{equation*}
and the wished convergence follows from the dominated convergence theorem.

Let us now study the case in which $t<s\searrow t$. For this end we notice
that, for all $s\geq t$, 
\begin{align}
& E\left[ \left\vert \left( e^{\left\vert t-s\right\vert A^{\ast }}-I\right)
\Phi _{t\vee s}\right\vert ^{2}\right]  \notag \\
& \leq c\left( E\left[ \left\vert \left( e^{\left( s-t\right) A^{\ast
}}-I\right) e^{\left( T-s\right) A^{\ast }}\xi \right\vert ^{2}\right]
\right.  \notag \\
& +E\left[ \left\vert \int_{s}^{T}\left( e^{\left( s-t\right) A^{\ast
}}-I\right) \left( C_{1}e^{\left( r-s\right) A}\right) ^{\ast
}Z_{r}dr\right\vert ^{2}\right]  \notag \\
& +E\left[ \left\vert \int_{s}^{T}\left( e^{\left( s-t\right) A^{\ast
}}-I\right) e^{\left( r-s\right) A^{\ast }}C_{2}^{\ast }Z_{r}dr\right\vert
^{2}\right]  \notag \\
& \left. +E\left[ \left\vert \int_{s}^{T}\left( e^{\left( s-t\right) A^{\ast
}}-I\right) e^{\left( r-s\right) A^{\ast }}Z_{r}dW_{r}\right\vert ^{2}\right]
\right)  \notag \\
& =I_{1}\left( s\right) +I_{2}\left( s\right) +I_{3}\left( s\right)
+I_{4}\left( s\right) .  \label{rl}
\end{align}%
For the first term we get from the dominated convergence theorem that%
\begin{equation*}
I_{1}\left( s\right) \leq kE\left[ \left\vert \left( e^{\left( s-t\right)
A^{\ast }}-I\right) \xi \right\vert ^{2}\right] \rightarrow 0\text{ as }%
s\searrow t.
\end{equation*}%
Next, 
\begin{align}
& I_{2}\left( s\right) \leq \left( \int_{s}^{T}\left( r-s\right) ^{-2\gamma
}dr\right) \times  \notag \\
\qquad & \times E\int_{t}^{T}I_{\left] s,T\right] }(r)\left\vert \left(
r-s\right) ^{\gamma }\left( e^{\left( s-t\right) A^{\ast }}-I\right) \left(
C_{1}e^{\left( r-s\right) A}\right) ^{\ast }Z_{r}\right\vert ^{2}dr.
\label{i2}
\end{align}%
We let $t<r\leq T$ and choose an arbitrary $s_{0}\in \left] t,r\right[ .$
Then, for all $t<s<s_{0},$ 
\begin{align*}
& \left\vert \left( e^{\left( s-t\right) A^{\ast }}-I\right) \left(
C_{1}e^{\left( r-s\right) A}\right) ^{\ast }Z_{r}\right\vert \\
& =\left\vert \left( e^{\left( s-t\right) A^{\ast }}-I\right) e^{\left(
s_{0}-s\right) A^{\ast }}\left( \left( C_{1}e^{\left( r-s_{0}\right)
A}\right) ^{\ast }Z_{r}\right) \right\vert \\
& \leq k\left\vert \left( e^{\left( s-t\right) A^{\ast }}-I\right) \left(
C_{1}e^{\left( r-s_{0}\right) A}\right) ^{\ast }Z_{r}\right\vert .
\end{align*}%
Obviously, the latter expression converges to $0$ as $s\searrow t.$
Consequently%
\begin{equation*}
I_{\left] s,T\right] }(r)\left\vert \left( r-s\right) ^{\gamma }\left(
e^{\left( s-t\right) A^{\ast }}-I\right) e^{\left( r-s\right) A^{\ast
}}C_{2}^{\ast }Z_{r}\right\vert ^{2}\underset{s\searrow t}{\rightarrow }0,%
\text{ for all }r>t,
\end{equation*}%
and from the dominated convergence theorem it follows that 
\begin{equation*}
I_{2}(s)\rightarrow 0\text{ as }s\searrow t.
\end{equation*}%
A similar argument yields $I_{3}(s)\rightarrow 0\text{ as }s\searrow t.$
Finally, for the last term, we have%
\begin{align*}
I_{4}(s)& \leq E\left[ \int_{s}^{T}\left\vert \left( e^{\left( s-t\right)
A^{\ast }}-I\right) e^{\left( r-s\right) A^{\ast }}Z_{r}\right\vert ^{2}dr%
\right] \\
& \leq E\left[ \int_{s}^{T}\left\vert \left( e^{\left( s-t\right) A^{\ast
}}-I\right) Z_{r}\right\vert ^{2}dr\right] ,
\end{align*}%
and, again by the dominated convergence theorem,%
\begin{equation*}
I_{4}(s)\rightarrow 0\text{ as }s\searrow t.
\end{equation*}%
Therefore, returning to (\ref{rl}) we get%
\begin{equation*}
\lim_{s\searrow t}E\left[ \left\vert \left( e^{\left\vert t-s\right\vert
A^{\ast }}-I\right) \Phi _{t\vee s}\right\vert ^{2}\right] =0.
\end{equation*}%
This concludes the proof of our lemma.$_{\blacksquare }$
\end{proof}

After having studied the existence and unique for the BSDE adjoint to our
forward stochastic control problem we are able now to characterize their
duality.

For the sake of simplicity, we shall assume from now on that $C_{1}$ takes
its values in $L_{2}\left( \Xi;H\right) .$

\begin{proposition}
Let $X^{x,u}$ be the unique mild solution of (\ref{eq1}) associated to an
admissible control $u$, and let $\left( Y,Z\right) $ be the unique mild
solution of (\ref{eq3}). Then the following duality relation holds true%
\begin{equation}
E\left[ \left\langle X_{T}^{x,u},Y_{T}\right\rangle \right] =E\left[
\left\langle x,Y_{0}\right\rangle \right] +E\left[ \int_{0}^{T}\left\langle
Bu_{s},Y_{s}\right\rangle ds\right] .  \label{eq0.4}
\end{equation}
\end{proposition}

\begin{proof}
For the proof of the duality relation we have the same difficulty as in the
proof of Theorem 1: we can't apply It\^{o}'s formula directly to our forward
SDE and our BSDE in infinite dimensions. This is why we consider the
following approximating equations%
\begin{equation*}
\left\{ 
\begin{array}{l}
dX_{t}^{n,\delta }=\left( A_{n}X_{t}^{n,\delta }+Bu_{t}\right) dt+\left(
C_{1}e^{\delta A}J_{n}+C_{2}\right) X_{t}^{n,\delta }dW_{t}, \\ 
X_{0}^{n}=x\in H,%
\end{array}%
\right.
\end{equation*}%
and%
\begin{equation*}
\left\{ 
\begin{array}{l}
dY_{t}^{n,\delta }=-\left( A_{n}^{\ast }Y_{t}^{n,\delta }+J_{n}^{\ast
}e^{\delta A^{\ast }}C_{1}^{\ast }Z_{t}^{n,\delta }+C_{2}^{\ast
}Z_{t}^{n,\delta }\right) dt+Z_{t}^{n,\delta }dW_{t}, \\ 
Y_{T}^{n,\delta }=\xi \in L^{2}\left( \Omega ,\mathcal{F}_{T},P;H\right) .%
\end{array}%
\right.
\end{equation*}%
Recall that $A_{n}^{\ast }=n(nI-A^{\ast })^{-1}A^{\ast }=J_{n}^{\ast
}A^{\ast }.$ To the above approximating equations we now can apply It\^{o}'s
formula, and we get%
\begin{equation}
E\left\langle Y_{s}^{n,\delta },X_{s}^{n,\delta }\right\rangle
=E\left\langle Y_{t}^{n,\delta },X_{t}^{n,\delta }\right\rangle +E\left[
\int_{t}^{s}\left\langle Bu_{r},Y_{r}^{n,\delta }\right\rangle dr\right] ,
\label{ineq0.5}
\end{equation}%
for all $0\leq t<s\leq T.$ Moreover, standard SDE and BSDE estimates show
that there exists some positive constant $k$ (not depending on $\delta $ and 
$n$), such that 
\begin{align*}
& E\left[ \sup_{t\in \lbrack 0,T]}\left\vert X_{t}^{n,\delta }\right\vert
^{2}\right] \leq k\left( 1+\left\vert x\right\vert ^{2}\right) \text{ and }
\\
& \sup_{t\in \lbrack 0,T]}E\left[ \left\vert Y_{t}^{n,\delta }\right\vert
^{2}\right] +E\left[ \int_{0}^{T}\left\vert Z_{s}^{n,\delta }\right\vert
^{2}ds\right] \leq kE\left[ \left\vert \xi \right\vert ^{2}\right] \text{.}
\end{align*}%
It follows that there exists some subsequence, still denoted $\left(
X^{n,\delta },Y^{n,\delta },Z^{n,\delta }\right) ,$ which converges weakly
to some limit $\left( X^{^{\prime }},Y^{^{\prime }},Z\right) $ in

$L^{2}\left( \Omega \times \lbrack 0,T];P\otimes dt;H\times H\right) $ $%
\times L^{2}\left( \Omega \times \lbrack 0,T];P\otimes dt;L_{2}\left( \Xi
;H\right) \right) $ as $n\rightarrow \infty ,$ $\delta \searrow 0.$ We
denote by $X$ the continuous version of $X^{^{\prime }}$; it is the unique
mild solution of equation (\ref{eq1}). Moreover, we let $Y$ be the $dtdP$%
-version of $Y^{\prime }$, which belongs to $C\left( \left[ 0,T\right]
;L^{2}\left( \Xi ;H\right) \right) $, and is, together with the process $Z$,
the unique mild solution of (\ref{eq3}). Moreover, from the above estimates
satisfied by $\left( X^{n,\delta },Y^{n,\delta },Z^{n,\delta }\right) $ we
get with the help of Mazur's theorem estimate (\ref{ineq0.4}) and 
\begin{equation*}
E\left[ \sup_{t\in \lbrack 0,T]}\left\vert X_{t}\right\vert ^{2}\right] \leq
k\left( 1+\left\vert x\right\vert ^{2}\right) \text{.}
\end{equation*}%
Moreover, if we take the weak limit as $n\rightarrow \infty $ and $\delta
\searrow 0$ in (\ref{ineq0.5}) we get 
\begin{equation*}
E\left\langle Y_{s}^{\prime },X_{s}^{\prime }\right\rangle =E\left\langle
Y_{t}^{\prime },X_{t}^{\prime }\right\rangle +E\left[ \int_{t}^{s}\left%
\langle Bu_{r},Y_{r}^{\prime }\right\rangle dr\right] ,\text{ }dtds\text{%
-a.e.},\text{ }0\leq t<s\leq T.
\end{equation*}%
Consequently,%
\begin{equation*}
E\left\langle Y_{s},X_{s}\right\rangle =E\left\langle
Y_{t},X_{t}\right\rangle +E\left[ \int_{t}^{s}\left\langle
Bu_{r},Y_{r}\right\rangle dr\right] ,\text{ for all }0\leq t<s\leq T.
\end{equation*}%
Finally, by taking $s=T$ and $t=0$, we get the assertion. The proof is
complete.$_{\blacksquare }$
\end{proof}

\medskip

The connection between equation (\ref{eq3}) and the approximate
controllability of (\ref{eq1}) is given by the following result that
generalizes those of the finite dimensional case.

\begin{proposition}
(i) The linear stochastic equation (\ref{eq1}) is approximately controllable
if and only if, for every finite time horizon $T>0,$ any solution of the
dual equation (\ref{eq3}) that satisfies $B^{\ast}Y_{s}=0$ $dP$-a.s., for
all $0\leq s\leq T$, necessarily vanishes $dsdP-a.s.,$ i.e. $Y_{s}=0$ $dP$%
-a.s., for all $0\leq s\leq T.$

(ii) The linear stochastic equation (\ref{eq1}) is approximately
null-controllable if and only if, for all finite time horizon $T>0,$ any
solution of the dual equation (\ref{eq3}) satisfying $B^{\ast}Y_{s}=0$ $dP$%
-a.s., for all $0\leq s\leq T,$ is such that $Y_{0}=0$ $dP$-a.e.
\end{proposition}

\begin{proof}
For any arbitrarily fixed time horizon $T>0$ we get from the previous
proposition that%
\begin{equation}
E\left[ \left\langle X_{T}^{x,u},Y_{T}\right\rangle \right] =E\left[
\left\langle x,Y_{0}\right\rangle \right] +E\left[ \int_{0}^{T}\left\langle
Bu_{s},Y_{s}\right\rangle ds\right] .  \label{eq9}
\end{equation}
We introduce the linear operator $M:\mathcal{U}\longrightarrow L^{2}\left(
\Omega,\mathcal{F}_{T},P;H\right) $ which associates to every admissible
control $u$ the mild solution of (\ref{eq1}) starting from $x=0$: 
\begin{equation*}
M(u)=X_{T}^{0,u}=\int_{0}^{T}e^{sA}Bu_{s}ds+%
\int_{0}^{T}e^{sA}CX_{s}^{0,u}dW_{s}.
\end{equation*}
Obviously, the approximate controllability (at time $T$) for (\ref{eq1}) is
equivalent to the condition that $M$ has an image space dense in $%
L^{2}(\Omega,\mathcal{F}_{T},P;H)$. This allows to deduce from (\ref{eq9})
the form of the dual operator of $M,$%
\begin{equation*}
M^{\ast}\xi=B^{\ast}Y.
\end{equation*}
On the other hand, since the density of the value domain of the bounded
linear operator $M\in L\left( L^{2}(\Omega,\mathcal{F}_{T},P;H)\right) $ is
equivalent with the condition that the kernel of its adjoint operator $%
M^{\ast}$ is trivial, we obtain from the above relation the first assertion.

For the proof of the second assertion we introduce the operator $%
L:H\longrightarrow$ $L^{2}\left( \Omega,\mathcal{F},P;H\right) $ which
associates to each initial state $x\in H$ the mild solution of (\ref{eq1})
corresponding to the control $u\equiv0:$%
\begin{equation*}
L(x)=e^{tA}x+\int_{0}^{T}e^{sA}CX_{s}^{x,0}dW_{s}.
\end{equation*}
From the relation $X_{T}^{x,u}=L(x)+M(u)$ we deduce easily that the
approximate null-controllability of $X$ is equivalent to the condition that $%
\overline{\mbox{Im}}(L)\subset\overline{\mbox{Im}}(M)$ ($\overline {\mbox{Im}%
}(L),\overline{\mbox{Im}}(M)$ are the closures of the image spaces of $L$
and $M$, resp.) and hence also to the following condition: 
\begin{equation*}
Ker\left( M^{\ast}\right) \subset Ker(L^{\ast}).
\end{equation*}
On the other hand, from (\ref{eq9}) we get $L^{\ast}\xi=Y_{0}$. This
relation together with $M^{\ast}\xi=B^{\ast}Y=0$ allow now to see the
equivalence between the approximate null-controllability of $X$ and the
condition given in the second assertion.$_{\blacksquare}$
\end{proof}

In what follows we will need the notion of the backward viability kernel
introduced by Buckdahn, Quincampoix, R\u{a}\c{s}canu \cite{bqr}

\begin{definition}
Let $K$ be a nonempty, convex, closed subset of $H.$

(i) A continuous stochastic process $\left\{ Y_{t},\text{ }t\in\left[ 0,T%
\right] \right\} $ is called viable in $K$ if and only if $Y_{t}\in K,$ $P$%
-a.s., for all $t\in\lbrack0,T].$

(ii) We say that the set $K$ enjoys the backward stochastic viability
property at time $T $ with respect to (\ref{eq3}) if for every $K$-valued
terminal condition $\eta\in L^{2}\left( \Omega,\mathcal{F}_{T},P;K\right) ,$
the solution $\left\{ Y_{t}^{\eta},\text{ }t\in\left[ 0,T\right] \right\} $
of (\ref{eq3}) is viable in $K.$

(iii) The largest closed, convex subset of $K$ enjoying the backward
stochastic viability property is called the backward stochastic viability
kernel of $K$.
\end{definition}

The notion of the stochastic viability kernel allows to reformulate the
criterion for the approximate controllability, stated in Proposition 2:

\begin{proposition}
The linear stochastic equation (\ref{eq1}) is approximately controllable if
and only if, for every finite time horizon $T>0,$ the backward stochastic
viability kernel of $Ker$ $B^{\ast}=\left\{ y\in H:\text{ }B^{\ast
}y=0\right\} $ at time $T$ \ with respect to (\ref{eq3}) is the trivial
subspace $\left\{ 0\right\} .$
\end{proposition}

\begin{remark}
\textit{In the finite dimensional case, the backward equation (\ref{eq3})
may be interpreted as a forward controlled equation. Therefore, instead of
studying the backward viability kernel, one may choose to investigate
approximate controllability with the help of the (forward) viability kernel.
Riccati methods are well adapted to control problems and allow nice
characterizations of the (forward) viability kernel. The authors of \cite%
{bqt} use these methods and show that approximate controllability of (\ref%
{eq1}) is equivalent to the following invariance condition:}%
\begin{equation*}
\text{The largest }\left( A^{\ast};C^{\ast}\right) \text{-strictly invariant
linear subspace of }Ker\text{ }B^{\ast}\text{ is }\left\{ 0\right\} .
\end{equation*}
\textit{We recall that a linear subspace }$V\subset%
\mathbb{R}
^{n}$\textit{\ is said to be }$(A^{\ast};C^{\ast})$\textit{-strictly
invariant if }$A^{\ast}V\subset Span\{V;C^{\ast}V\}=\left\{ \lambda v+\mu
w:v\in V,\text{ }w\in C^{\ast}V\right\} $.

\textit{If }$H$\textit{\ is infinite dimensional, and }$A$\textit{\ is a
generator of a strongly continuous group, similar arguments apply.}
\end{remark}

\begin{remark}
\textit{Let us suppose that the Brownian motion }$W$\textit{\ is
1-dimensional, }$B\in \mathcal{L}(H)$\textit{, and }$C$\textit{\ is a linear
(possibly unbounded) operator on }$H$\textit{\ such that }$A^{\ast }B^{\ast
}=B^{\ast }A^{\ast }$\textit{\ and }$B^{\ast }C^{\ast }=C^{\ast }B^{\ast }.$%
\textit{\ Then (\ref{eq1}) is approximately controllable if and only if the
image space }$\func{Im}(B)$\textit{\ is dense in }$H$\textit{.}

\textit{Indeed, let us notice that if }$(Y,Z)$\textit{\ is the mild solution
of (\ref{eq3}) and satisfies (\ref{ineq0.4}), then}%
\begin{equation*}
Y_{t}=e^{(T-t)A^{\ast }}\xi +\int_{t}^{T}e^{(s-t)A^{\ast }}C^{\ast
}Z_{s}ds-\int_{t}^{T}e^{(s-t)A^{\ast }}Z_{s}dW_{s},
\end{equation*}%
\textit{and, from the commutativity of }$B^{\ast \text{ }}$\textit{with }$%
A^{\ast }$\textit{\ and with }$C^{\ast },$%
\begin{equation*}
B^{\ast }Y_{t}=e^{(T-t)A^{\ast }}B^{\ast }\xi +\int_{t}^{T}e^{(s-t)A^{\ast
}}C^{\ast }B^{\ast }Z_{s}ds-\int_{t}^{T}e^{(s-t)A^{\ast }}B^{\ast
}Z_{s}dW_{s}.
\end{equation*}%
\textit{Thus, }$B^{\ast }Y_{t}$\textit{\ is the unique mild solution of the
following BSDE:}%
\begin{equation*}
\left\{ 
\begin{array}{l}
d\widetilde{Y}_{t}=-A^{\ast }\widetilde{Y}_{t}dt-C^{\ast }\widetilde{Z}%
_{t}dt+\widetilde{Z}_{t}dW_{t}, \\ 
\widetilde{Y}_{T}=B^{\ast }\xi .%
\end{array}%
\right. 
\end{equation*}%
\textit{Obviously, }$\widetilde{Y}=0$\textit{\ if and only if }$B^{\ast }\xi
=0$\textit{\ }$P$\textit{-a.s.}$.$\textit{\ Thus, from Proposition 2 it
follows that Eq. (\ref{eq1}) is approximately controllable if, for all }$\xi
\in L^{2}\left( \Omega ,\mathcal{F}_{T},P;H\right) ,$\textit{\ the relation }%
$B^{\ast }\xi =0,$\textit{\ }$P-a.s.,$\textit{\ implies that }$\xi =0,$%
\textit{\ }$P-a.s.$\textit{\ This is, of course, equivalent with the density
of the image space }$\func{Im}(B)$\textit{\ in }$\mathit{H}$\textit{.}
\end{remark}

\section{A necessary condition for approximate controllability}

We have seen that approximate controllability for the forward controlled
equation (\ref{eq1}) is equivalent to the following (approximate)
observability condition on the dual equation (\ref{eq3}) :%
\begin{equation}
"B^{\ast}Y_{t}=0,\text{ }dP-a.s.,\text{ for all }t\in\left[ 0,T\right] ,%
\text{ implies }Y_{T} =0,\text{ }dP-a.s."  \label{suff}
\end{equation}
In the deterministic case, Russell and Weiss \cite{rw} generalized the
Hautus test of observability for infinite dimensional equations with an
operator $A $ that is supposed to generate an exponentially stable
semigroup. In what follows we assume besides (A1) and (A2) the following
additional condition:

\textbf{(A3)}\textit{\ The linear operator }$A$\textit{\ generates an
exponentially stable, strongly continuous semigroup of operators.}

Under the assumptions (A1)-(A3) we can prove the following statement:

\begin{proposition}
A necessary condition for the approximate controllability of (\ref{eq1}) is
that, for every $y\in D\left( A^{\ast}\right) $ and every $\alpha<0,$ 
\begin{equation}
\left\vert B^{\ast}y\right\vert +\left\vert \left( A^{\ast}-\alpha I\right)
y\right\vert >0,\text{ whenever }y\neq0.  \tag{N1}  \label{n1}
\end{equation}
\end{proposition}

\begin{proof}
In order to prove the claim, let us first notice that $H_{1}=D\left(
A\right) $ endowed with the norm $\left\vert h\right\vert _{1}=\left\vert
\left( A^{\ast }-\alpha I\right) h\right\vert _{H}$ is a Hilbert space. It
is well known that, under the above assumptions, the family of norms indexed
by $\alpha <0$ are equivalent with the usual graph norm on $H_{1}.$ For
every $y\in D\left( A^{\ast }\right) $ we let $\left( Y^{y},Z^{y}\right) $
denote the unique mild solution in $H$ of the BSDE 
\begin{equation*}
\left\{ 
\begin{array}{l}
dY_{t}^{y}=-A^{\ast }Y_{t}^{y}dt-C^{\ast }Z_{t}^{y}dt+Z_{t}^{y}dW_{t}, \\ 
Y_{T}=y.%
\end{array}%
\right.
\end{equation*}%
Since all data of this BSDE is deterministic it is immediate that $Y^{y}$ is
deterministic and $Z^{y}=0$. In particular, we see that $Y_{t}^{y}=e^{\left(
T-t\right) A^{\ast }}y$ is a classical solution (in $H$) of%
\begin{equation*}
\left\{ 
\begin{array}{l}
dY_{t}^{y}=-\alpha Y_{t}^{y}dt-e^{\left( T-t\right) A^{\ast }}\left( A^{\ast
}-\alpha I\right) ydt, \\ 
Y_{T}^{\eta }=y,%
\end{array}%
\right.
\end{equation*}%
and the function $B^{\ast }Y^{y}$ is a classical solution of the following
equation:%
\begin{equation*}
\left\{ 
\begin{array}{l}
d\left( B^{\ast }Y_{t}^{y}\right) =-\alpha \left( B^{\ast }Y_{t}^{y}\right)
dt-B^{\ast }e^{\left( T-t\right) A^{\ast }}\left( A^{\ast }-\alpha I\right)
ydt \\ 
B^{\ast }Y_{T}^{y}=B^{\ast }y.%
\end{array}%
\right.
\end{equation*}%
It follows easily from this equation that $B^{\ast }Y_{t}^{y}=0,$ for all $%
t\in \left[ 0,T\right] ,$ if and only if%
\begin{equation*}
\left\{ 
\begin{array}{l}
B^{\ast }y=0, \\ 
B^{\ast }e^{tA^{\ast }}\left( A^{\ast }-\alpha I\right) y=0,\text{ for all }%
t\in \left[ 0,T\right] .%
\end{array}%
\right.
\end{equation*}%
Consequently, the condition (\ref{suff}) gives the following necessary
condition for the approximate controllability of (4): 
\begin{equation*}
"B^{\ast }Y_{t}^{y}=0,\text{ for all }t\in \left[ 0,T\right] ,\text{ implies 
}y=0.\text{ }"
\end{equation*}%
Obviously, the two latter conditions allow to conclude that%
\begin{equation}
\left\{ 
\begin{array}{l}
B^{\ast }y=0,\text{ } \\ 
B^{\ast }e^{tA^{\ast }}\left( A^{\ast }-\alpha I\right) y=0,\text{ for all }%
t\in \left[ 0,T\right] ,%
\end{array}%
\right. \text{implies }y=0\text{,}  \label{nec}
\end{equation}%
and the estimate 
\begin{equation*}
\left\vert B^{\ast }e^{tA^{\ast }}\left( A^{\ast }-\alpha I\right)
y\right\vert \leq k\left\vert \left( A^{\ast }-\alpha I\right) y\right\vert ,
\end{equation*}%
in combination with (\ref{nec}) allows to complete the proof.$_{\blacksquare
}$
\end{proof}

\begin{remark}
\textit{Jacob, Partington \cite{jp} studied the approximate controllability
for a deterministic system. They supposed}

\textit{(JP) }$A$\textit{\ is an infinitesimal generator of an exponentially
stable, strongly continuous semigroup which possesses a sequence of
normalized eigenvectors }$\left\{ e_{i}\right\} $\textit{\ corresponding to
the eigenvalues }$\left\{ \lambda_{i}\right\} $\textit{\ such that }$\sup
_{i}\lambda_{i}<0.$\textit{\ \ Moreover, they considered the case of a
1-dimensional input space, i.e. }$B\in L\left( 
\mathbb{R}
;H\right) $\textit{.}

\textit{In this particular case, the necessary and sufficient condition for
approximate controllability of the deterministic system}%
\begin{equation*}
\left\{ 
\begin{array}{l}
dX_{t}^{x,u}=\left( AX_{t}^{x,u}+Bu_{t}\right) dt, \\ 
X_{0}=x\in H,%
\end{array}
\right.
\end{equation*}
\textit{found by the authors, says that for all }$y\in H_{1}$\textit{\ and
all }$\alpha<0,$\textit{\ }%
\begin{equation*}
\left\vert B^{\ast}y\right\vert ^{2}+\left\vert \left( A^{\ast}-\alpha
I\right) y\right\vert ^{2}>0\text{ whenever }y\neq0.
\end{equation*}
\end{remark}

\medskip

\begin{remark}
\textit{For the case in which }$H$\textit{\ is $n$-dimensional Euclidean
space (stochastic) approximate controllability was studied by Buckdahn,
Quincampoix, Tessitore \cite{bqt} and Goreac \cite{g}. The equivalent
condition for approximate controllability reads }%
\begin{equation}
\text{The largest }\left( A^{\ast};C^{\ast}\right) \text{-strictly invariant
subspace of }Ker\text{ }B^{\ast}\text{ is }\left\{ 0\right\} .  \label{bqt}
\end{equation}

\textit{Let us suppose that, for the framework studied by these authors,
there exists a bounded linear operator }$D\in L(U)$\textit{\ such that }$%
B^{\ast }C^{\ast}=DB^{\ast}$\textit{. Then we get that }$Ker$\textit{\ }$%
B^{\ast}$\textit{\ is }$C^{\ast}$\textit{- invariant, and thus (\ref{bqt})
can be written as follows:}%
\begin{equation}
\text{The largest }A^{\ast}\text{-invariant subspace of \ }Ker\text{ }%
B^{\ast }\text{ is }\left\{ 0\right\} .  \label{bqtpart}
\end{equation}

\textit{\ Moreover, under the assumptions of Jacob, Partington \cite{jp}
(JP), it is obvious that (\ref{n1}) is equivalent to (\ref{bqtpart}).
Indeed, if (\ref{n1}) holds true, then } 
\begin{equation*}
\left\{ 
\begin{array}{l}
B^{\ast}e_{i}\neq0,\text{ for all }1\leq i\leq n, \\ 
\lambda_{i}\neq\lambda_{j},\text{ for all }1\leq i,j\leq n,\text{ }i\neq j.%
\end{array}
\right.
\end{equation*}
\textit{(see Jacob, Partington \cite{jp}, Theorem 4.1). Let }$V$\textit{\
denote the largest }$A^{\ast}$\textit{-invariant subspace of \ }$Ker$\textit{%
\ }$B^{\ast}$\textit{, and let us suppose that there exists some linear
combination } $v=\sum_{k=1}^{m}v_{i_{k}}e_{i_{k}}$ \textit{\ such that } $%
v\in V$, \textit{\ where } $m\leq n,i_{k}\in\left\{ 1,2,\ldots ,n\right\} $ 
\textit{\ and } $v_{i_{k}}\neq0$, \textit{\ for all }$1\leq k\leq m.$ 
\textit{\ Then, for all} $j\geq
m-1,\,\sum_{k=1}^{m}\lambda_{i_{k}}^{j}v_{i_{k}}e_{i_{k}}\in V$. \textit{\
Thus, since }%
\begin{equation*}
\det\left[ \lambda_{i_{k}}^{j}v_{i_{k}}\right] _{k,j}=\prod\limits_{1\leq
k\leq m}v_{i_{k}}\prod\limits_{1\leq k<j\leq m}\left(
\lambda_{i_{j}}-\lambda_{i_{k}}\right) \neq0,
\end{equation*}
\textit{we get that }%
\begin{equation*}
span\left\{ e_{i_{k}},1\leq k\leq m\right\} \subset V.
\end{equation*}
\textit{It follows that }$V=span\left\{ e_{i_{k}},1\leq k\leq N\right\} ,$%
\textit{\ for some }$N\leq n.$\textit{\ But then }$B^{\ast}e_{i_{k}}=0,$%
\textit{\ and this contradicts our assumption and we have that }$V=\left\{
0\right\} .$

\textit{For the converse, if (\ref{bqtpart}) holds true and }$y\in H_{1}$%
\textit{\ such that }%
\begin{equation*}
\left\vert B^{\ast}y\right\vert ^{2}+\left\vert \left( A^{\ast}-\alpha
I\right) y\right\vert ^{2}=0\text{ \ for some }\alpha<0,
\end{equation*}
\textit{then }$V=span\left\{ y\right\} $\textit{\ is }$A^{\ast}$\textit{%
-invariant and included in }$Ker$\textit{\ }$B^{\ast}.$\textit{\ It follows
that }$y=0,$\textit{\ and we get (\ref{n1}). This latter argument applies
also when }$H$\textit{\ has infinite dimension.}
\end{remark}

Let us now make the following assumptions:

\textbf{(B) }\textit{\ }$W$\textit{\ is supposed to be a 1-dimensional
Brownian motion, the control state space }$U$\textit{\ is a bounded closed
subspace of some separable real Hilbert space }$V$\textit{, }$B\in \mathcal{L%
}(V;H),$\textit{\ }$A$\textit{\ is a self adjoint operator which generates a
semigroup of contractions on }$H$\textit{, and the operator }$C$\textit{\
admits a decomposition } 
\begin{equation*}
C=C_{1}+C_{2},
\end{equation*}
\textit{\ of two linear operators } $C_{1},\,C_{2}$ \textit{\ which are
supposed to have the following properties:}

\textit{1) }$C_{2}$\textit{\ is a bounded operator from }$H$\textit{\ to }$H$%
\textit{;}

\textit{2) for all }$t>0,$\textit{\ }$C_{1}e^{tA}$\textit{, }$e^{tA}C_{1}\in%
\mathcal{L}\left( H\right) .$\textit{\ Moreover, we suppose that there exist
some }$\gamma\in\left[ 0,\frac{1}{2}\right) $\textit{\ and some positive
constant }$L>0$\textit{\ such that }%
\begin{equation*}
\left\vert C_{1}e^{tA}\right\vert _{\mathcal{L}\left( H\right) }+\left\vert
e^{tA}C_{1}\right\vert _{\mathcal{L}\left( H\right) }\leq Lt^{-\gamma},
\end{equation*}
\textit{for all }$t>0.$

\textit{3) There exists some constant }$a>\frac{1}{2}$\textit{\ such that }%
\begin{equation*}
\text{ }A+aC_{1}^{\ast}C_{1}\text{ is dissipative}.
\end{equation*}

We recall the following

\begin{definition}
Let $A$ be the generator of a $C_{0}-$semigroup on the Hilbert space $H$ and 
$C$ is a linear operator on $H$. We say that $C$ is a \textit{class-}$%
\mathcal{P}$ \textit{perturbation }of $A$ if $C$ is closed, 
\begin{equation*}
D\left( C\right) \supset\cup_{t>0}e^{tA}\left( H\right) \mbox{ and
} \int_{0}^{1}\left\vert Ce^{tA}\right\vert dt<\infty.
\end{equation*}
\end{definition}

Obviously, under the above assumptions, the operator $C$ is a class-$%
\mathcal{P}$ perturbation of $A$. It follows that $A+\lambda C$ is the
generator of a $C_{0}$-semigroup $\left( e^{t\left( A+\lambda C\right)
}\right) _{t\geq0}$ for all $\lambda\in%
\mathbb{R}
$ (cf. Davies \cite{d} Theorem 3.5)$.$

\medskip

For the study of the main result of this section we will need the following
estimates:

\begin{lemma}
Under our standard assumptions we have that, for some constant $k$, 
\begin{equation*}
\left\vert C_{1}e^{t\left( A+\lambda C\right) }\right\vert _{\mathcal{L}%
(H)}+\left\vert e^{t\left( A+\lambda C\right) }C_{1}\right\vert _{\mathcal{L}%
(H)}\leq k\left( t^{-\gamma}+1\right) ,
\end{equation*}
for all $t\in\left[ 0,T\right] .$
\end{lemma}

\begin{proof}
From the theory of general perturbation of generators it follows that%
\begin{align*}
e^{t\left( A+\lambda C\right) }x & =e^{tA}x+\lambda\int_{0}^{t}e^{\left(
t-s\right) A}C_{1}e^{s\left( A+\lambda C\right) }x \\
& +\lambda\int_{0}^{t}e^{\left( t-s\right) A}C_{2}e^{s\left( A+\lambda
C\right) }x,
\end{align*}
for all $x\in H$. Then, by applying on both sides of the above relation the
bounded operator $C_{2}$, we get the following norm estimate: 
\begin{align*}
\left\vert C_{1}e^{t\left( A+\lambda C\right) }x\right\vert & \leq
t^{-\gamma}\left\vert x\right\vert +\lambda\int_{0}^{t}\left( t-s\right)
^{-\gamma}\left\vert C_{1}e^{s\left( A+\lambda C\right) }x\right\vert ds \\
& +k\int_{0}^{t}\left( t-s\right) ^{-\gamma}\left\vert x\right\vert ds,
\end{align*}
for all $t\in\lbrack0,T]$. Here $k$ denotes again a generic constant which
can depend on $\lambda$ and $T.$ Thus, applying Cauchy-Schwarz inequality
yields 
\begin{align*}
\left\vert C_{1}e^{t\left( A+\lambda C\right) }x\right\vert ^{2} & \leq
k\left( \left( t^{-2\gamma}+t^{2-2\gamma}\right) \left\vert x\right\vert
^{2}+t^{1-2\gamma}\int_{0}^{t}\left\vert C_{1}e^{s\left( A+\lambda C\right)
}x\right\vert ^{2}ds\right) \\
& \leq k\left( \left( t^{-2\gamma}+1\right) \left\vert x\right\vert
^{2}+\int_{0}^{t}\left\vert C_{1}e^{s\left( A+\lambda C\right) }x\right\vert
^{2}ds\right) ,
\end{align*}
and from Gronwall's inequality we finally get%
\begin{equation*}
\left\vert C_{1}e^{t\left( A+\lambda C\right) }x\right\vert ^{2}\leq k\left(
t^{-\gamma}+1\right) ^{2}\left\vert x\right\vert ^{2}.
\end{equation*}
It follows that $C_{1}e^{t\left( A+\lambda C\right) }\in\mathcal{L}\left(
H\right) $ and 
\begin{equation*}
\left\vert C_{1}e^{t\left( A+\lambda C\right) }\right\vert _{\mathcal{L}%
(H)}\leq k\left( t^{-\gamma}+1\right) ,
\end{equation*}
for all $t\in\left[ 0,T\right] .$ Using a similar argument we can prove that 
$e^{t\left( A+\lambda C\right) }C_{1}\in\mathcal{L}\left( H\right) $ and 
\begin{equation*}
\left\vert e^{t\left( A+\lambda C\right) }C_{1}\right\vert _{\mathcal{L}%
(H)}\leq k\left( t^{-\gamma}+1\right) ,
\end{equation*}
for all $t\in\left[ 0,T\right] ._{\blacksquare}$
\end{proof}

\medskip

To establish the main result of this section we shall further introduce the
following set standing for the joint dissipativity condition on $A,C$: 
\begin{equation*}
\Lambda =\left\{ \lambda \in 
\mathbb{R}
:\exists a>\frac{1}{2}\text{ such that }A+\lambda C_{1}+aC_{1}^{\ast }C_{1}%
\text{ is dissipative}\right\} .
\end{equation*}

\begin{remark}
1. If $C\in\mathcal{L}\left( H\right) $ is a bounded operator, then $\Lambda=%
\mathbb{R}
.$

2. $\Lambda$ contains at least the origin $\left\{ 0\right\} $.

3. If $C_{1}$ is dissipative and the assumption (B) holds true, then $%
\mathbb{R}
_{+}\subset\Lambda$.
\end{remark}

We now can state our main result of this section.

\begin{theorem}
Under assumption \textit{(B)}, a necessary condition for the approximate
controllability of (\ref{eq1}) is%
\begin{equation}
\left\vert B^{\ast}y\right\vert +\left\vert \left( A^{\ast}+\lambda C^{\ast
}-\alpha I\right) y\right\vert >0\text{, for all }y\neq0\text{, and all }%
\left( \lambda,\alpha\right) \in\Lambda\times%
\mathbb{R}
_{-}.  \label{N2}
\end{equation}
\end{theorem}

The above necessary condition is an immediate consequence of Proposition 4
and a $\lambda$-wise application of the following result:

\begin{theorem}
If (\ref{eq1}) is approximately controllable, then the system%
\begin{equation}
\left\{ 
\begin{array}{l}
dX_{t}=\left( \left( A+\lambda C\right) X_{t}+Bv_{t}\right) dt+\left(
C+\lambda I\right) X_{t}dW_{t}, \\ 
X_{0}=x\in H,%
\end{array}%
\right.   \label{eqC}
\end{equation}%
which is governed by the control process $v\in L_{\tciFourier }^{2}\left( %
\left[ 0,T\right] ;V\right) $ is also approximately controllable.
\end{theorem}

\begin{proof}
\underline{Step 1}. Approximation of (\ref{eqC}) by an equation with bounded
operators admitting the application of It\^{o}'s formula.

\medskip

For all $u\in L_{\tciFourier}^{0}\left( \left[ 0,T\right] ;U\right) $, we
denote by $X_{n,\delta }^{x,u}$ the unique mild solution of the controlled
forward equation 
\begin{equation*}
\left\{ 
\begin{array}{l}
dX_{n,\delta}^{x,u}(t)=A_{n}X_{n,\delta}^{x,u}(t)dt+Bu\left( t\right)
dt+J_{n}^{\ast}e^{\delta A^{\ast}}Ce^{\delta
A}J_{n}X_{n,\delta}^{x,u}(t)dW_{t}, \\ 
X_{n,\delta}^{x,u}(0)=x\in H,%
\end{array}
\right.
\end{equation*}
where $J_{n}=$ $\left( I-n^{-1}A\right) ^{-1}$ and $A_{n}=J_{n}A$. This
approximation of the operators $A$ (by $A_{n}$) and $C$ (by $J_{n}^{\ast
}e^{\delta A^{\ast}}Ce^{\delta A}J_{n}$) explains by the same difficulties
we have already met in the proof of Theorem 1. Our special choice of the
approximation allows to conserve the joint dissipativity condition also for
the approximating operators and allows now to apply It\^{o}'s formula.

Let $\mathcal{E}\left( \lambda W\right) $ denote the Dol\'ean-Dade
exponential of $\lambda W$, i.e., $\mathcal{E}\left( \lambda W\right)
_{t}:=e^{\lambda W_{t}-\frac{\lambda^{2}}{2}t},\, t\in\left[ 0,T\right] $.
Then, from It\^{o}'s formula applied to $\mathcal{E}\left( \lambda W\right)
_{t}X_{n,\delta}^{x,u}(t)$ it follows that%
\begin{equation*}
\left\{ 
\begin{array}{l}
d\left( \mathcal{E}\left( \lambda W\right) _{t}X_{n,\delta}^{x,u}(t)\right)
=\left( A_{n}+\lambda J_{n}^{\ast}e^{\delta A^{\ast}}Ce^{\delta
A}J_{n}\right) \left( \mathcal{E}\left( \lambda W\right) _{t}X_{n,\delta
}^{x,u}(t)\right) dt \\ 
\text{ \ \ \ \ \ \ \ \ \ \ \ \ \ \ \ \ \ \ \ \ \ }+B\left( \mathcal{E}\left(
\lambda W\right) _{t}u\left( t\right) \right) dt \\ 
\text{ \ \ \ \ \ \ \ \ \ \ \ \ \ \ \ \ \ \ \ \ \ }+\left( J_{n}^{\ast
}e^{\delta A^{\ast}}Ce^{\delta A}J_{n}+\lambda I\right) \left( \mathcal{E}%
\left( \lambda W\right) _{t}X_{n,\delta}^{x,u}(t)\right) dW_{t}, \\ 
X_{n,\delta}^{x,u}(0)=x\in H.%
\end{array}
\right.
\end{equation*}

After the above application of It\^{o}'s formula we would like to take the
limit as $n\rightarrow+\infty$ and then as $\delta\downarrow0$ in order to
get an equation which coincides with that we would get if we applied
formally It\^{o}'s formula to $\mathcal{E}\left( \lambda W\right)
_{t}X^{x,u}(t),$ where $X^{x,u}$ denotes the unique mild solution of (4).
For taking these limits we need the following result whose proof will be
given later.

\medskip

\begin{proposition}
Under the assumptions on Theorem 2 and with the notations introduced above
we have that, for all $x\in H$, 
\begin{equation}
\lim_{n}\sup_{0\leq t\leq T}\left\vert e^{t\left( A_{n}+\lambda J_{n}^{\ast
}e^{\delta A^{\ast}}Ce^{\delta A}J_{n}\right) }x-e^{t\left( A+\lambda
e^{\delta A^{\ast}}Ce^{\delta A}\right) }x\right\vert =0,\, \delta >0,
\label{rel1}
\end{equation}
and%
\begin{equation}
\lim_{\delta}\sup_{0\leq t\leq T}\left\vert e^{t\left( A+\lambda e^{\delta
A^{\ast}}Ce^{\delta A}\right) }x-e^{t\left( A+\lambda C\right) }x\right\vert
=0.  \label{rel2}
\end{equation}
\end{proposition}
\end{proof}

We continue the

\begin{proof}
of our theorem. With the help of the above proposition we are now able to
prove

\medskip

\underline{Step 2}. Let $X^{x,u}$ denote the unique mild solution of (4).
Then the process $\mathcal{E} \left( \lambda W\right) _{\cdot}X^{x,u}\left(
\cdot\right) $ is the unique mild solution of (\ref{eqC}). Moreover, 
\begin{equation}
\sup_{0\leq t\leq T}E\left[ \left\vert \mathcal{E}\left( \lambda W\right)
_{t}X^{x,u}\left( t\right) \right\vert ^{p}\right] \leq c_{p}\left(
1+\left\vert x\right\vert ^{p}\right) .  \label{estim}
\end{equation}

\medskip

For proving this statement we first notice that from standard estimates, for
all $p>2$,%
\begin{align*}
E\left[ \sup_{0\leq t\leq T}\left\vert \mathcal{E}\left( \lambda W\right)
_{t}X_{n,\delta }^{x,u}(t)\right\vert ^{p}\right] & \leq c_{p}\left(
1+\left\vert x\right\vert ^{p}\right) ,\text{ and} \\
E\left[ \sup_{0\leq t\leq T}\left\vert X_{n,\delta }^{x,u}(t)\right\vert ^{p}%
\right] & \leq c_{p}\left( 1+\left\vert x\right\vert ^{p}\right) ;
\end{align*}%
$c_{p}$ denotes a generic constant independent of $n,$ $\delta $ and $u\in
L_{\tciFourier }^{0}\left( \left[ 0,T\right] ;U\right) $. Then, for any $%
\delta >0,$ there exists a subsequence of

$\left( \mathcal{E}\left( \lambda W\right) _{\cdot }X_{n,\delta
}^{x,u}(\cdot ),X_{n,\delta }^{x,u}(\cdot )\right) _{n}$ , still denoted by $%
\left( \mathcal{E}\left( \lambda W\right) _{\cdot }X_{n,\delta }^{x,u}(\cdot
),X_{n,\delta }^{x,u}(\cdot )\right) _{n}$ ,which converges in the weak
topology on

$L^{p}\left( \left[ 0,T\right] \times \Omega ;H\right) \times L^{2p}\left( %
\left[ 0,T\right] \times \Omega ;H\right) $ to some limit $\left( X_{\delta
}^{\prime }\left( \cdot \right) ,X_{\delta }^{\prime \prime }\left( \cdot
\right) \right) $. With the help of Proposition 5 we can show that $%
X_{\delta }^{\prime }$ is a unique mild solution of 
\begin{equation*}
\left\{ 
\begin{array}{l}
dX_{\delta }^{\prime }\left( t\right) =\left( A+\lambda e^{\delta A^{\ast
}}Ce^{\delta A}\right) X_{\delta }^{\prime }(t)dt \\ 
\text{ \ \ \ \ \ \ \ \ \ \ }+B\left( \mathcal{E}\left( \lambda W\right)
_{t}u\left( t\right) \right) dt+\left( e^{\delta A^{\ast }}Ce^{\delta
A}+\lambda I\right) X_{\delta }^{\prime }(t)dW_{t}, \\ 
X_{\delta }^{^{\prime }}(0)=x\in H,%
\end{array}%
\right.
\end{equation*}%
and $X_{\delta }^{\prime \prime }$ is a mild solution of 
\begin{equation}
\left\{ 
\begin{array}{l}
dX_{\delta }^{\prime \prime }(t)=\left( AX_{\delta }^{\prime \prime
}(t)+Bu_{t}\right) dt+e^{\delta A^{\ast }}Ce^{\delta A}X_{\delta }^{\prime
\prime }(t)dW_{t}, \\ 
X_{\delta }^{\prime \prime }(0)=x\in H.%
\end{array}%
\right.
\end{equation}%
On the other hand, it follows from the general theory of SDEs in infinite
dimensions that these mild solutions are unique and that 
\begin{equation}
\sup_{0\leq t\leq T}E\left[ \left\vert X_{\delta }^{\prime }\left( t\right)
\right\vert ^{p}\right] \leq c_{p}\left( 1+\left\vert x\right\vert
^{p}\right) .  \label{estp}
\end{equation}%
Moreover, taking into account that

$\mathcal{E}\left( \lambda W\right) _{\cdot}\zeta\left( \cdot\right) \in L^{%
\frac{2p}{2p-1}}\left( \left[ 0,T\right] \times\Omega;H\right) $, for all $%
\zeta\in L^{\frac{p}{p-1}}\left( \left[ 0,T\right] \times \Omega;H\right) $,

\noindent we get 
\begin{align*}
& E\left[ \int_{0}^{T}\left\langle X_{\delta}^{\prime}(t),\zeta\left(
t\right) \right\rangle dt\right] =\lim_{n}E\left[ \int_{0}^{T}\left\langle 
\mathcal{E}\left( \lambda W\right) _{t}X_{n,\delta}^{x,u}(t),\zeta\left(
t\right) \right\rangle dt\right] \\
& =\lim_{n}E\left[ \int_{0}^{T}\left\langle X_{n,\delta}^{x,u}(t),\mathcal{E}%
\left( \lambda W\right) _{t}\zeta\left( t\right) \right\rangle dt\right] \\
& =E\left[ \int_{0}^{T}\left\langle \mathcal{E}\left( \lambda W\right)
_{t}X_{\delta}^{\prime\prime}(t),\zeta\left( t\right) \right\rangle dt\right]
.
\end{align*}

This relation allows to identify the processes $X_{\delta}^{\prime}\left(
\cdot\right) $ and $\mathcal{E}\left( \lambda W_{\cdot}\right) X_{\delta
}^{\prime\prime}\left( \cdot\right) $ as elements of $L^{p}\left( \left[ 0,T%
\right] \times\Omega;H\right) $. Moreover, if $X_{\delta}^{x,u}$ denotes the
continuous version of $X_{\delta}^{\prime\prime}$ and $\widetilde {X}%
_{\delta}^{x,u}$the continuous version of $X_{\delta}^{\prime}$, we have%
\begin{equation*}
\widetilde{X}_{\delta}^{x,u}\left( t\right) =\mathcal{E}\left( \lambda
W\right) _{t}X_{\delta}^{x,u}\left( t\right) ,\text{ }dP\text{-a.s, for all }%
t\in\left[ 0,T\right] ,
\end{equation*}
and inequality (\ref{estp}) takes the form 
\begin{equation*}
\sup_{0\leq t\leq T}E\left[ \left\vert \mathcal{E}\left( \lambda W\right)
_{t}X_{\delta}^{x,u}\left( t\right) \right\vert ^{p}\right] \leq c_{p}\left(
1+\left\vert x\right\vert ^{p}\right) .
\end{equation*}
By repeating the argument for letting $\delta\rightarrow0$ we get the result
stated in step 2.

\smallskip

After having related equation (4) with equation (\ref{eqC}) we can prove now
the theorem in its proper sense.

\medskip

\underline{Step 3}. Conclusion.

\medskip

If $\xi\in L^{2}\left( \Omega,\mathcal{F}_{T},P;H\right) $, then, for every $%
\varepsilon>0$ there exists some $\xi^{\varepsilon}\in L^{\infty}\left(
\Omega,\mathcal{F}_{T},P;H\right) $ such that 
\begin{equation*}
E\left[ \left\vert \xi^{\varepsilon}-\xi\right\vert ^{2}\right]
\leq\varepsilon.
\end{equation*}
It follows from (\ref{estim}) that the family

$\left\{ \left\vert \mathcal{E}\left( \lambda W\right) _{T}X^{x,u}\left(
T\right) -\xi^{\varepsilon}\right\vert ^{2},\text{ }u\in
L_{\tciFourier}^{0}\left( \left[ 0,T\right] ;U\right) \right\} $

\noindent is uniformly integrable. Consequently, there exists $%
M_{\varepsilon }>0$ such that 
\begin{equation*}
E\left[ \left\vert \mathcal{E}\left( \lambda W\right) _{T}X^{x,u}\left(
T\right) -\xi^{\varepsilon}\right\vert ^{2}1_{\left\{ \mathcal{E}\left(
\lambda W\right) _{T}>M_{\varepsilon}\right\} }\right] \leq\varepsilon,
\end{equation*}
for all $u\in L_{\tciFourier}^{0}\left( \left[ 0,T\right] ;U\right) $. If
the equation (\ref{eq1}) is approximately controllable, then there exists $%
u_{\varepsilon}\in L_{\tciFourier}^{0}\left( \left[ 0,T\right] ;U\right) $
such that 
\begin{equation*}
E\left[ \left\vert X_{T}^{x,u_{\varepsilon}}-\xi^{\varepsilon}\mathcal{E}%
\left( \lambda W\right) _{T}^{-1}\right\vert ^{2}\right] \leq \frac{%
\varepsilon}{M_{\varepsilon}^{2}},
\end{equation*}
and we get%
\begin{align*}
E\left[ \left\vert \mathcal{E}\left( \lambda W\right)
_{T}X_{T}^{x,u_{\varepsilon}}-\xi^{\varepsilon}\right\vert ^{2}\right] &
\leq M_{\varepsilon}^{2}E\left[ \left\vert
X_{T}^{x,u_{\varepsilon}}-\xi^{\varepsilon}\mathcal{E}\left( \lambda
W\right) _{T}^{-1}\right\vert ^{2}\right] \\
& +E\left[ \left\vert \mathcal{E}\left( \lambda W\right)
_{T}X^{x,u_{\varepsilon}}\left( T\right) -\xi^{\varepsilon}\right\vert
^{2}1_{\left\{ \mathcal{E}\left( \lambda W\right) _{T}>M_{\varepsilon
}\right\} }\right] \\
& \leq2\varepsilon.
\end{align*}
Therefore, also (\ref{eqC}) is approximately controllable. The proof of our
theorem is now complete.$_{\blacksquare}$
\end{proof}

\medskip

However, the proof of Proposition 5 still remains open:

\medskip

\begin{proof}
(of {Proposition 5}). Due to the definition of the approximation of the
operators $A$ and $C$ given in step 1 of the proof of the above theorem we
have for all $x\in\mathcal{D}\left( A+\lambda e^{\delta A^{\ast}}Ce^{\delta
A}\right) $, 
\begin{equation}
\lim_{n}\left( A_{n}+\lambda J_{n}^{\ast}e^{\delta A^{\ast}}Ce^{\delta
A}J_{n}\right) x=\left( A+\lambda e^{\delta A^{\ast}}Ce^{\delta A}\right) x.
\label{2}
\end{equation}
For all $n,$ the operator $A_{n}+\lambda J_{n}^{\ast}e^{\delta
A^{\ast}}Ce^{\delta A}J_{n}$ is bounded. Therefore, it generates a $C_{0}$%
-semigroup $\left( e^{t\left( A_{n}+\lambda J_{n}^{\ast}e^{\delta
A^{\ast}}Ce^{\delta A}J_{n}\right) }\right) _{t}$ and the application $%
t\longmapsto\left\vert e^{t\left( A_{n}+\lambda J_{n}^{\ast}e^{\delta
A^{\ast}}Ce^{\delta A}J_{n}\right) }\right\vert $ is continuous. From the
general theory of perturbation of generators, we have%
\begin{align*}
e^{t\left( A_{n}+\lambda J_{n}^{\ast}Ce^{\delta A}J_{n}\right) }x &
=e^{tA_{n}}x \\
& +\lambda\int_{0}^{t}e^{\left( t-s\right) A_{n}}J_{n}^{\ast}e^{\delta
A^{\ast}}C_{1}e^{\delta A}J_{n}e^{s\left( A_{n}+\lambda
J_{n}^{\ast}e^{\delta A^{\ast}}Ce^{\delta A}J_{n}\right) }xds \\
& +\lambda\int_{0}^{t}e^{\left( t-s\right) A_{n}}J_{n}^{\ast}e^{\delta
A^{\ast}}C_{2}e^{\delta A}J_{n}e^{s\left( A_{n}+\lambda
J_{n}^{\ast}e^{\delta A^{\ast}}Ce^{\delta A}J_{n}\right) }xds.
\end{align*}
It follows that, for $n$ great enough 
\begin{equation*}
\left\vert e^{t\left( A_{n}+\lambda J_{n}^{\ast}Ce^{\delta A}J_{n}\right)
}\right\vert \leq1+\lambda\int_{0}^{t}\left( \delta^{-\gamma}+k\right)
\left\vert e^{s\left( A_{n}+\lambda J_{n}^{\ast}e^{\delta
A^{\ast}}Ce^{\delta A}J_{n}\right) }\right\vert ds,
\end{equation*}
where $k>0$ is a generic constant (which may depend on $\delta$ but not on $%
n $), and Gronwall's inequality yields%
\begin{equation}
\left\vert e^{t\left( A_{n}+\lambda J_{n}^{\ast}e^{\delta
A^{\ast}}Ce^{\delta A}J_{n}\right) }\right\vert \leq e^{kt},  \label{3}
\end{equation}
for all $t>0$, and all $n\in%
\mathbb{N}
$. Then, from (\ref{2}) and (\ref{3}) we get (cf. Davies \cite{d} Th. 3.17)
that (\ref{rel1}) holds true, for all $\delta>0$ and all $x\in\mathcal{D}(A)$%
.

To prove the second assertion, we notice that%
\begin{align*}
e^{t\left( A+\lambda e^{\delta A^{\ast}}Ce^{\delta A}\right) }x &
=e^{tA}x+\int_{0}^{t}e^{\left( t-s\right) A}\lambda e^{\delta
A^{\ast}}C_{1}e^{\delta A}e^{s\left( A+\lambda e^{\delta A^{\ast}}Ce^{\delta
A}\right) }xds \\
& +\int_{0}^{t}e^{\left( t-s\right) A}\lambda e^{\delta
A^{\ast}}C_{2}e^{\delta A}e^{s\left( A+\lambda e^{\delta A^{\ast}}Ce^{\delta
A}\right) }xds,
\end{align*}
for all $x\in H$. Then, recalling that $A$ is self adjoint,we obtain%
\begin{align*}
\left\vert e^{t\left( A+\lambda e^{\delta A^{\ast}}Ce^{\delta A}\right)
}x\right\vert & \leq\left\vert x\right\vert +\lambda\int_{0}^{t}\left\vert
e^{s\left( A+\lambda e^{\delta A^{\ast}}Ce^{\delta A}\right) }x\right\vert
\left( t-s\right) ^{-\gamma}ds \\
& +k\int_{0}^{t}\left\vert e^{s\left( A+\lambda e^{\delta
A^{\ast}}Ce^{\delta A}\right) }x\right\vert ds
\end{align*}
($k$ is again a generic constant independent of $\delta$). Thus, with the
notation%
\begin{equation*}
f\left( t\right) =\left\vert e^{t\left( A+\lambda e^{\delta
A^{\ast}}Ce^{\delta A}\right) }x\right\vert ,
\end{equation*}
the latter estimate takes the form%
\begin{equation*}
f\left( t\right) \leq\left\vert x\right\vert +\lambda\int_{0}^{t}f(s)\left(
t-s\right) ^{-\gamma}ds+k\int_{0}^{t}f(s)ds.
\end{equation*}
Then, by Cauchy-Schwarz inequality,%
\begin{align}
f\left( t\right) & \leq\left\vert x\right\vert +k\left( \frac {t^{\frac{%
1-2\gamma}{2}}}{\sqrt{1-2\gamma}}+t^{\frac{1}{2}}\right) \left(
\int_{0}^{t}f^{2}(s)ds\right) ^{\frac{1}{2}}  \label{rel0} \\
& \leq\left\vert x\right\vert +k\left( T^{\frac{1}{2}}\vee1\right) \left(
\int_{0}^{t}f^{2}(s)ds\right) ^{\frac{1}{2}},
\end{align}
and, consequently, 
\begin{equation*}
f^{2}(t)\leq2\left( \left\vert x\right\vert ^{2}+k\left( T\vee1\right)
\int_{0}^{t}f^{2}(s)ds\right) .
\end{equation*}
To the latter estimate we apply Gronwall's inequality and take the square
root after. This yields%
\begin{equation*}
f(t)\leq\sqrt{2}\left\vert x\right\vert e^{k\left( T\vee1\right) t}.
\end{equation*}
Therefore, from the definition of $f(t)$ it follows that 
\begin{equation}
\sup_{\delta>0}\left\vert e^{t\left( A+\lambda e^{\delta A^{\ast}}Ce^{\delta
A}\right) }\right\vert \leq Me^{ct},  \label{rel3}
\end{equation}
for all $t\leq T$, where $M$ and $c$ are positive constants that are
independent of $\delta>0$. On the other hand, for all $x\in\mathcal{D}\left(
A+\lambda C\right) $ we have 
\begin{equation}
\lim_{\delta\rightarrow0}\left( A+\lambda e^{\delta A^{\ast}}Ce^{\delta
A}\right) x=\left( A+\lambda C\right) x.  \label{rel4}
\end{equation}
The second assertion follows (cf. Davies \cite{d} Th. 3.17).$_{\blacksquare}$
\end{proof}

\medskip

In the following we discuss tow examples to illustrate the results of this
section.

\medskip

\begin{example}
Given a regular domain $\mathcal{O\subset%
\mathbb{R}
}^{N}$ we consider the following stochastic partial differential equation%
\begin{equation}
\left\{ 
\begin{array}{c}
d_{t}X^{u}\left( t,x\right) =\sum_{i,j=1}^{N}\partial_{i}\left(
a_{i,j}(x)\partial_{j}X^{u}\left( t,x\right) \right) dt+u\left( t\right)
b(x)dt \\ 
+\sum_{i=1}^{N}c_{i}(x)\partial_{i}X^{u}(t,x)dW_{t}, \\ 
\multicolumn{1}{l}{X^{u}(t,x)=0,\text{ }\forall\left( t,x\right) \in\left[
0,T\right] \times\partial\mathcal{O}\text{,}} \\ 
\multicolumn{1}{l}{X^{u}(0,x)=\xi(x),\text{ }\forall x\in\mathcal{O}\text{,}}%
\end{array}
\right.  \label{SPDE}
\end{equation}
where $u$ is an admissible control process taking its values in $%
\mathbb{R}
$ . We suppose that $a(x)\left( =(a_{i,j}(x)\right) \sigma(x)\sigma^{\ast
}(x)$ for some $C_{\ell,b}^{\infty}$ matrix $\sigma$ of $N\times N$-type, $%
c=\left( c_{1},\ldots,c_{N}\right) \in C_{\ell,b}^{\infty}\left( \mathcal{O};%
\mathbb{R}
^{N}\right) ,$ $b\in H^{1}\left( \mathcal{O}\right) $ and $\xi\in
L^{2}\left( \Omega,\mathcal{F}_{T},P;L^{2}\left( \mathcal{O}\right) \right)
. $ Moreover, we suppose that the couple of coefficients $(a,c)$ satisfies
the standard ellipticity condition%
\begin{equation}
\sum_{i,j=1}^{N}\left( a_{i,j}(x)-\alpha c_{i}(x)c_{j}(x)\right) \lambda
_{i}\lambda_{j}\geq0,\text{ }  \label{el}
\end{equation}
for some $\alpha>\frac{1}{2}$ and for all $\lambda\in%
\mathbb{R}
^{N}.$ Then, if we put 
\begin{align*}
& H=L^{2}\left( \mathcal{O}\right) , \\
& \mathcal{D}(A)=H^{2}\left( \mathcal{O}\right) \cap H_{0}^{1}\left( 
\mathcal{O}\right) ,\text{ }A\zeta=\sum_{i,j=1}^{N}\partial_{i}\left(
a_{i,j}(x)\partial_{j}\zeta\left( x\right) \right) , \\
& \mathcal{D}\left( C\right) =H^{1}\left( \mathcal{O}\right) ,\text{ }%
C\zeta=c\cdot\nabla\zeta,
\end{align*}
we get that 
\begin{equation*}
\mathcal{D}\left( C^{\ast}\right) =H^{1}\left( \mathcal{O}\right) ,\text{ }%
C^{\ast}\zeta=-c\cdot\nabla\zeta-\zeta\sum_{i=1}^{N}\partial_{i}c_{i}.
\end{equation*}
The ellipticity condition (\ref{el}) insures that the dual backward
stochastic partial differential equation%
\begin{equation}
\left\{ 
\begin{array}{c}
d_{t}Y(t,x)=-\left( \sum_{i,j=1}^{N}\partial_{i}\left(
a_{i,j}(x)\partial_{j}Y\left( t,x\right) \right) \right) dt+\left( \sum
_{i=1}^{N}c_{i}(x)\partial_{i}Z(t,x)\right) dt \\ 
+\left( \sum_{i=1}^{N}\partial_{i}c_{i}(x)Z(t,x)\right) dt+Z\left(
t,x\right) dW_{t}, \\ 
\multicolumn{1}{l}{Y(t,x)=Z(t,x)=0,\text{ }\forall\left( t,x\right) \in\left[
0,T\right] \times\partial\mathcal{O}\text{,}} \\ 
\multicolumn{1}{l}{Y(T,x)=\eta(x),\text{ }\forall x\in\mathcal{O}\text{,}}%
\end{array}
\right.  \label{BSPDE}
\end{equation}
has a unique mild solution. Thus we know that the approximate
controllability of (\ref{SPDE}) is equivalent to the approximate
observability of (\ref{BSPDE}).

From (\ref{n1}) it follows that, if (\ref{SPDE}) is approximately
controllable and if $\zeta_{n}(x)$ is a complete orthonormal base consisting
of eigenvectors for $A,$ then every coefficient of $b$ in this base must be
non null.
\end{example}

\begin{remark}
The problem of controllability for the deterministic version of (\ref{SPDE})
has been treated by Carleman estimates method in Fursikov, Imanuvilov \cite%
{fi}.
\end{remark}

The condition (N2) is non trivially more general then (N1) as proven by the
following

\begin{example}
We consider the following equation%
\begin{equation}
\left\{ 
\begin{array}{l}
d_{t}X^{u}\left( t,x\right) =\triangle X^{u}\left( t,x\right) dt+u\left(
t\right) b\left( x\right) dt \\ 
\text{ \ }+\left( 2\sin \left( \pi x\right) \int_{0}^{1}X^{u}\left(
t,y\right) \sin \left( \pi y\right) dy\right) dW_{t}, \\ 
X^{u}\left( t,0\right) =X^{u}\left( t,1\right) =0,\text{ }\forall t\in \left[
0,T\right] . \\ 
X^{u}\left( 0,x\right) =\xi \left( x\right) ,\text{ }\forall x\in \left(
0,1\right) ,%
\end{array}%
\right.   \label{EX2}
\end{equation}%
where $u$ is an admissible real-valued bounded control process and $b\in
L^{2}\left( 0,1\right) $. This equation can be expressed as an infinite
dimensional linear equation. For this we put%
\begin{align*}
& H=L^{2}\left( 0,1\right) ,\text{ }\mathcal{D}(A)=H^{2}\left( 0,1\right)
\cap H_{0}^{1}\left( 0,1\right) , \\
& A\zeta =\triangle \zeta ,\text{ for all }\zeta \in \mathcal{D}(A), \\
& C\zeta \left( \cdot \right) =2\sin \left( \pi \cdot \right)
\int_{0}^{1}\zeta \left( y\right) \sin \left( \pi y\right) dy,\text{ for all 
}\zeta \in H.
\end{align*}%
Obviously $C$ is a self-adjoint bounded linear operator on $H$. Furthermore,
suppose that 
\begin{equation*}
b_{n}=\sqrt{2}\int_{0}^{1}b\left( y\right) \sin \left( \pi y\right) dy\neq 0,
\end{equation*}%
for all $n\geq 1$. Then (N1) is obviously satisfied. However, if we choose $%
\lambda =-3\pi ^{2}$, $\alpha =-4\pi ^{2}$ and $\zeta \left( \cdot \right) =-%
\frac{b_{2}\sqrt{2}}{b_{1}}\sin \left( \pi \cdot \right) +\sqrt{2}\sin
\left( 2\pi \cdot \right) $, we have%
\begin{equation*}
\left\vert \left( A^{\ast }+\lambda C^{\ast }-\alpha I\right) \zeta
\right\vert ^{2}+\left\vert B^{\ast }\zeta \right\vert ^{2}=0.
\end{equation*}%
It follows that (N2) is not satisfied which implies that the equation (\ref%
{EX2}) cannot be approximately controllable.
\end{example}


\begin{thebibliography}{99}
\bibitem{brt} V. Barbu, A. R\u{a}\c{s}canu, G. Tessitore (2003), \textit{%
Carleman estimates and Controllability of stochastic heat equations with
multiplicative noise}, Appl. Math. Optim. 47:97--120, pp.98-120.

\bibitem{bqr} Buckdahn, R., Quincampoix, M., Rascanu, A. (2000), \textit{%
Viability property for a backward stochastic differential equation and
applications to partial differential equations}, Probab. Theory Relat.
Fields 116, No.4, pp. 485-504.

\bibitem{bqt} Buckdahn, R., Quincampoix, M., Tessitore, G. (2006), \textit{A
Characterization of Approximately Controllable Linear Stochastic
Differential Equations, } Stochastic Partial Differential Equations and
Applications, G. Da Prato and L. Tubaro Eds\ Series of Lecture Notes in pure
and appl. Math., Chapman \& Hall Vol.245, pp. 253-260.

\bibitem{c} Confortola, F. (2004), \textit{Dissipative backward stochastic
differential equations in infinite dimensions}.

\bibitem{dpz} Da Prato, G., Zabczyk, J. (1992), \textit{Stochastic equations
in infinite dimensions}, Cambridge University Press, Cambridge.

\bibitem{dpz2} Da Prato, G., Zabczyk, J. (1996), \textit{Ergodicity for
infinite-dimensional systems. London Mathematical Society Lecture Note
Series, 229}, Cambridge University Press, Cambridge.

\bibitem{d} Davies, E. B. (1980), \textit{One-parameter semigroups}, London
Mathematical Society Monographs, 15. Academic Press, Inc., London-New York.

\bibitem{fcgar} Fernandez-Cara, E., Garrido-Atienza, M. J. Real J. (1999), 
\textit{On the approximate controllability of a stochastic parabolic
equation with a multiplicative noise, }C. R. Acad. Sci. Paris, t. 328, Serie
I, pp. 675-680.

\bibitem{ft} Fuhrman, M., Tessitore, G. (2002), \textit{Nonlinear Kolmogorov
equations in infinite dimesional spaces: the backward stochastic
differential equations approach and applications to optimal control}, Ann.
Probab. 30, pp. 1397-1465.

\bibitem{fi} Fursikov, A., Imanuvilov, O. (1996), \textit{Controllability of
evolution equations}, vol. 34. Seoul National University.

\bibitem{g} Goreac, D. (2007), \textit{Approximate Controllability for
Linear Stochastic Differential Equations with Control Acting on the Noise, }%
Applied Analysis and Differential Equations, Ia\c{s}i, Rom\^{a}nia 4 - 9
September 2006, World Scientific Publishing, pp. 153-164.

\bibitem{hmy} Hu, Y., Ma, J., Yong, J. (2002), \textit{On semi-linear
degenerate backward stochastic partial differential equations}, Probability
Theory and Related Fields, vol. 123, no. 3, pp. 381--411.

\bibitem{jp} Jacob, B., Partington, J., R. (2006), \textit{On
controllability of diagonal systems with one-dimensional input space, }%
Systems and Control Letters 55, pp. 321 -- 328.

\bibitem{jz} Jacob, B., Zwart, H. (2001), \textit{Exact observability of
diagonal systems with a finitedimensional output operator}, Systems Control
Lett. 43 101--109.

\bibitem{lp} Liu, Y., Peng, S. (2002), \textit{Infinite horizon backward
stochastic differential equation and exponential convergence index
assignment of stochastic control systems}, Automatica, 38, pp. 1417-1423.

\bibitem{my1} Ma, J., Yong, J. (1997), \textit{Adapted solution of a
degenerate backward SPDE, with applications}, Stochastic Processes and Their
Applications, vol. 70, no. 1, pp. 59--84.

\bibitem{my2} Ma, J., Yong, J. (1999), \textit{On linear, degenerate
backward stochastic partial differential equations}, Probability Theory and
Related Fields, vol. 113, no. 2, pp. 135--170.

\bibitem{pp} Pardoux, E., Peng, S.G. (1990), \textit{Adapted solutions of a
backward stochastic differential equation,} Systems and Control Letters, 14,
pp. 55-61.

\bibitem{p} Peng, S.G. (1994), \textit{Backward Stochastic Differential
Equation and Exact Controllability of Stochastic Control Systems}, Progr.
Natur. Sci. vol. 4, No. 3, pp. 274-284.

\bibitem{rw} Russell, D.L., Weiss, G. (1994), \textit{A general necessary
condition for exact observability}, SIAM J. Control Optim. 32 (1), pp. 1--23.

\bibitem{st} Sirbu, M, Tessitore, G. (2001), \textit{Null controllability of
an infinite dimensional SDE with state and control-dependent noise}, Systems
and Control Letters, 44, pp. 385-394.

\bibitem{t} Tessitore, G. (1996), \textit{Existence, uniqueness and space
regularity of the adapted solutions of a backward SPDE,} Stochastic Analysis
and Applications, vol. 14, no. 4, pp. 461--486.
\end{thebibliography}
\end{document}